\documentclass[a4paper]{amsart}
\usepackage{color}
\usepackage[latin1]{inputenc}
\usepackage[T1]{fontenc}
\usepackage{amsfonts}
\usepackage{amssymb}
\usepackage{amsmath}
\usepackage{amsthm}
\usepackage{stmaryrd}
\usepackage{enumerate}
\usepackage{multirow}
\usepackage{graphicx}

\usepackage{graphicx,type1cm,eso-pic,color}
\usepackage{pstricks}
\usepackage{float}
\theoremstyle{definition}
\newtheorem{theorem}{Theorem}
\newtheorem{lemma}{Lemma}

\newtheorem{assumption}{Assumption}
\newtheorem{remark}{Remark}

\begin{document}
\setlength\arraycolsep{2pt}
\title{On a Projection Least Squares Estimator for Jump Diffusion Processes}
\author{H\'el\`ene HALCONRUY$^{\dag,\diamond}$}
\email{helene.halconruy@devinci.fr}
\author{Nicolas MARIE$^{\diamond}$}
\email{nmarie@parisnanterre.fr}
\date{}
\maketitle
\noindent
$^{\dag}$L\'eonard de Vinci P\^ole universitaire, Research Center,\\
12 avenue L\'eonard de Vinci, 92400, Courbevoie, France.
\\
\\
$^{\diamond}$Laboratoire Modal'X, Universit\'e Paris Nanterre,\\
200 avenue de la R\'epublique, 92001, Nanterre, France.
\\
\\
{\bf Corresponding author:} Nicolas MARIE (nmarie@parisnanterre.fr).
%

% Abstract.

%
\begin{abstract}
This paper deals with a projection least squares estimator of the drift function of a jump diffusion process $X$ computed from multiple independent copies of $X$ observed on $[0,T]$. Risk bounds are established on this estimator and on an associated adaptive estimator. Finally, some numerical experiments are provided.
\end{abstract}
\tableofcontents
%

% Section : Introduction.

%
\section{Introduction}\label{introduction_section}
Let $Z = (Z_t)_{t\in [0,T]}$ be the compound Poisson process defined by
\begin{displaymath}
Z_t :=\sum_{n = 1}^{\nu_t}\zeta_n =
\int_{0}^{t}\int_{-\infty}^{\infty}z\mu(ds,dz)
\end{displaymath}
for every $t\in [0,T]$, where $\nu = (\nu_t)_{t\in [0,T]}$ is a (usual) Poisson process of intensity $\lambda > 0$, independent of the $\zeta_n$'s which are i.i.d. random variables of probability distribution $\pi$, and $\mu$ is the Poisson random measure of intensity $m(ds,dz) :=\lambda\pi(dz)ds$ defined by
\begin{displaymath}
\mu([0,t]\times dz) :=
|\{s\in [0,t] : Z_s - Z_{s^-}\in dz\}|
\textrm{ $;$ }
\forall t\in [0,T].
\end{displaymath}
In the sequel, $Z$ is {\it replaced} by the centered martingale $\mathfrak Z = (\mathfrak Z_t)_{t\in [0,T]}$ defined by
\begin{displaymath}
\mathfrak Z_t :=
Z_t -\int_{0}^{t}\int_{-\infty}^{\infty}zm(ds,dz) =
Z_t -\mathfrak c_{\zeta}\lambda t
\end{displaymath}
for every $t\in [0,T]$, where $\mathfrak c_{\zeta}$ is the (common) expectation of the $\zeta_n$'s. Now, let us consider the stochastic differential equation
\begin{equation}\label{main_equation}
X_t = x_0 +\int_{0}^{t}b(X_s)ds +
\int_{0}^{t}\sigma(X_s)dB_s +
\int_{0}^{t}\gamma(X_s)d\mathfrak Z_s
\textrm{ $;$ }
t\in [0,T],
\end{equation}
where $x_0\in\mathbb R$, $B = (B_t)_{t\in [0,T]}$ is a Brownian motion independent of $Z$, $b\in C^1(\mathbb R)$ and its derivative is bounded, and $\sigma,\gamma :\mathbb R\rightarrow\mathbb R$ are bounded Lipschitz continuous functions such that $\inf_{x\in\mathbb R}\sigma(x)^2\wedge\gamma(x)^2 > 0$. Under these conditions on $b$, $\sigma$ and $\gamma$, Equation (\ref{main_equation}) has a unique (strong) solution $X = (X_t)_{t\in [0,T]}$.

As for continuous diffusion processes, the major part of the estimators of the drift function in stochastic differential equations driven by jump processes are computed from one path of the solution to Equation (\ref{main_equation}) and converges when $T\rightarrow\infty$ (see Schmisser (2014), Gloter et al. (2018), Amorino et al. (2022), etc.). The existence and the uniqueness of the ergodic stationary solution to Equation (\ref{main_equation}) is then required, and obtained thanks to a restrictive dissipativity condition on $b$. For stochastic differential equations driven by a pure-jump L\'evy process, some authors have also studied estimation methods based on high frequency observations, on a fixed time interval, of one path of the solution (see Cl\'ement and Gloter (2019,2020)).

Now, consider $X^i :=\mathcal I(x_0,B^i,\mathfrak Z^i)$ for every $i\in\{1,\dots,N\}$, where $\mathcal I(.)$ is the It\^o map associated to Equation (\ref{main_equation}) and $(B^1,\mathfrak Z^1),\dots,(B^N,\mathfrak Z^N)$ are $N\in\mathbb N^*$ independent copies of $(B,\mathfrak Z)$. The estimation of the drift function $b$ from a continuous-time or a discrete-time observation of $(X^1,\dots,X^N)$ is a functional data analysis problem already investigated in the parametric and in the nonparametric frameworks for continuous diffusion processes (see Ditlevsen and De Gaetano (2005), Picchini and Ditlevsen (2011), Delattre et al. (2013), Comte and Genon-Catalot (2020b), Denis et al. (2021), Marie and Rosier (2023), etc.). Up to our knowledge, no such estimator of the drift function has been already proposed for jump diffusion processes. So, our paper deals with a projection least squares estimator $\widehat b_m$ of $b$ computed from $X^1,\dots,X^N$, which means that $\widehat b_m$ is minimizing the objective function
\begin{displaymath}
\tau\longmapsto
\gamma_N(\tau) :=
\frac{1}{NT}\sum_{i = 1}^{N}\left(
\int_{0}^{T}\tau(X_{s}^{i})^2ds - 2\int_{0}^{T}\tau(X_{s}^{i})dX_{s}^{i}
\right)
\end{displaymath}
on a $m$-dimensional function space $\mathcal S_m$. Precisely, risk bounds are established on $\widehat b_m$ and on the adaptive estimator $\widehat b_{\widehat m}$, where
\begin{displaymath}
\widehat m =\arg\min_{m\in\widehat{\mathcal M}_N}
\{-\|\widehat b_m\|_{N}^{2} + {\rm pen}(m)\}
\end{displaymath}
with $\widehat{\mathcal M}_N\subset\{1,\dots,N\}$,
\begin{displaymath}
{\rm pen}(m) :=\mathfrak c_{\rm cal}\frac{m}{N}
\textrm{ $;$ }
\forall m\in\mathbb N^*
\end{displaymath}
and $\mathfrak c_{\rm cal} > 0$ is a constant to calibrate in practice.

In Section \ref{section_projection_LS_estimator}, a detailed definition of the projection least squares estimator of $b$ is provided. Section \ref{section_risk_bound} deals with a risk bound on $\widehat b_m$ and Section \ref{section_model_selection} with a risk bound on the adaptive estimator $\widehat b_{\widehat m}$. Finally, some numerical experiments are provided in Section \ref{section_numerical_experiments}. The proofs (resp. tables and figures) are postponed to Appendix \ref{section_proofs} (resp. Appendix \ref{section_tables_figures}).
\\
\\
\textbf{Notations and basic definitions:}
\begin{itemize}
 \item $\mathfrak c_{\zeta^n} :=\mathbb E(\zeta_{1}^{n})$ for every $n\in\mathbb N^*$.
 \item Consider $d\in\mathbb N^*$. The $j$-th component of any $\mathbf x\in\mathbb R^d$ is denoted by $\mathbf x_j$ or $[\mathbf x]_j$.
 \item For every $k\in\mathbb N^*$, $\|.\|_{k,d}$ is the norm on $\mathbb R^d$ defined by
 \begin{displaymath}
 \|\mathbf x\|_{k,d} :=
 \left(\sum_{j = 1}^{d}|\mathbf x_j|^k\right)^{1/k}
 \textrm{ $;$ }
 \forall\mathbf x\in\mathbb R^d.
 \end{displaymath}
 \item The spectral norm on the space $\mathcal M_d(\mathbb R)$ of the $d\times d$ real matrices is denoted by $\|.\|_{\rm op}$:
 \begin{displaymath}
 \|{\bf A}\|_{\rm op} :=
 \sup_{x\in\mathbb R^d :\|x\|_{2,d} = 1}\|\mathbf Ax\|_{2,d}
 \textrm{ $;$ }
 \forall\mathbf A\in\mathcal M_d(\mathbb R).
 \end{displaymath}
\end{itemize}
%

% Section : A projection least squares estimator of the drift function.

%
\section{A projection least squares estimator of the drift function}\label{section_projection_LS_estimator}
%

% Subsection : The objective function.

%
\subsection{The objective function}\label{section_objective_function}
Assume that the probability distribution of $X_s$ has a density $p_s(x_0,.)$ with respect to Lebesgue's measure for every $s\in (0,T]$, that $s\mapsto p_s(x_0,x)$ belongs to $\mathbb L^1([0,T],dt)$ for every $x\in\mathbb R$ which legitimates to consider the density function $f_T$ defined by
\begin{displaymath}
f_T(x) :=\frac{1}{T}\int_{0}^{T}p_s(x_0,x)ds
\textrm{ $;$ }
\forall x\in\mathbb R,
\end{displaymath}
and that
\begin{displaymath}
\int_{-\infty}^{\infty}b(x)^4f_T(x)dx <\infty.
\end{displaymath}
%

% Remark : Sufficient condition for the solution to get a well-bounded density.

%
\begin{remark}\label{sufficient_condition_density}
Assume that $b$ and $\pi$ satisfy the following additional conditions:
\begin{enumerate}
 \item The function $b$ belongs to the Kato class
 \begin{displaymath}
 \mathbb K_2 :=
 \left\{\varphi :\mathbb R\rightarrow\mathbb R :
 \lim_{\delta\rightarrow 0}\sup_{x\in\mathbb R}
 \int_{0}^{\delta}\int_{-\infty}^{\infty}
 |\varphi(x + y) +\varphi(x - y)|s^{1/2}(|y| + s^{1/2})^{-3}dyds = 0
 \right\}.
 \end{displaymath}
 \item The L\'evy measure $\pi_{\lambda}(.) :=\lambda\pi(.)$ has a density $\theta$ with respect to Lebesgue's measure. Moreover, there exists $\alpha\in (0,2)$ such that $z\in\mathbb R\mapsto\theta(z)|z|^{1 +\alpha}$ is bounded, and if $\alpha = 1$, then
 \begin{displaymath}
 \int_{r < |z| < R}z\theta(z)dz = 0
 \textrm{ $;$ }
 \forall R > r > 0.
 \end{displaymath}
\end{enumerate}
By Chen et al. (2017), Theorem 1.1 and the remark p. 126, l. 5-7, in Amorino and Gloter (2021), for every $s\in (0,T]$, the probability distribution of $X_s$ has a density $p_s(x_0,.)$ with respect to Lebesgue's measure, and there exist two constants $\mathfrak c_p,\mathfrak m_p > 0$, not depending on $s$, such that
\begin{equation}\label{sufficient_condition_density_1}
p_s(x_0,x)\leqslant\mathfrak c_p\left[
s^{-1/2}\exp\left(-\mathfrak m_p\frac{(x - x_0)^2}{s}\right) +\frac{s}{(s^{1/2} +
|x - x_0|)^{1 +\alpha}}\right]
\textrm{ $;$ }
\forall x\in\mathbb R.
\end{equation}
So,
\begin{itemize}
 \item $f_T$ is well-defined and even bounded, which is crucial in Section \ref{section_model_selection}. Indeed, since $-1/2 < 1 - (1 +\alpha)/2 < 1/2$, for every $x\in\mathbb R$,
 \begin{displaymath}
 0\leqslant f_T(x)\leqslant
 \frac{\mathfrak c_p}{T}\left(\int_{0}^{T}s^{-1/2}ds +
 \int_{0}^{T}s^{1 - (1 +\alpha)/2}ds\right) =
 \frac{\mathfrak c_p}{T}\left[2T^{1/2} +\frac{T^{2 - (1 +\alpha)/2}}{2 - (1 +\alpha)/2}\right]
 <\infty.
 \end{displaymath}
 \item If there exists a constant $\mathfrak c_b > 0$ and $\varepsilon\in (0,\alpha)$ as close as possible to $0$ such that $|b(x)|\leqslant\mathfrak c_b(1 + |x|)^{(\alpha -\varepsilon)/4}$ for every $x\in\mathbb R$, then
 \begin{displaymath}
 b(x)^4p_s(x_0,x)\underset{x\rightarrow\pm\infty,s\rightarrow 0^+}{=}
 O\left(\frac{s^{-1/2}}{|x - x_0|^{1 +\varepsilon}}\right),
 \end{displaymath}
 which leads to
 \begin{displaymath}
 \int_{-\infty}^{\infty}b(x)^4f_T(x)dx <\infty.
 \end{displaymath}
\end{itemize}
\end{remark}
Now, let us consider the objective function $\gamma_N(.)$ defined by
\begin{displaymath}
\gamma_N(\tau) :=
\frac{1}{NT}\sum_{i = 1}^{N}\left(
\int_{0}^{T}\tau(X_{s}^{i})^2ds - 2\int_{0}^{T}\tau(X_{s}^{i})dX_{s}^{i}
\right)
\end{displaymath}
for every $\tau\in\mathcal S_m$, where $m\in\{1,\dots,N_T\}$, $N_T := [NT] + 1$, $\mathcal S_m := {\rm span}\{\varphi_1,\dots,\varphi_m\}$, $\varphi_1,\dots,\varphi_{N_T}$ are continuous functions from $I$ into $\mathbb R$ such that $(\varphi_1,\dots,\varphi_{N_T})$ is an orthonormal family in $\mathbb L^2(I,dx)$, and $I\subset\mathbb R$ is a non-empty interval. For any $\tau\in\mathcal S_m$,
\begin{eqnarray*}
 \mathbb E(\gamma_N(\tau)) & = &
 \frac{1}{T}\int_{0}^{T}\mathbb E(|\tau(X_s) - b(X_s)|^2)ds -
 \frac{1}{T}\int_{0}^{T}\mathbb E(b(X_s)^2)ds\\
 & = &
 \int_{-\infty}^{\infty}(\tau(x) - b(x))^2f_T(x)dx -
 \int_{-\infty}^{\infty}b(x)^2f_T(x)dx.
\end{eqnarray*}
Then, the closer $\tau$ is to $b$, the smaller $\mathbb E(\gamma_N(\tau))$. For this reason, the estimator of $b$ minimizing $\gamma_N(.)$ is studied in this paper.
%

% Subsection : The projection least squares estimator and some related matrices.

%
\subsection{The projection least squares estimator and some related matrices}
Consider
\begin{displaymath}
J :=
\sum_{j = 1}^{m}\theta_j\varphi_j
\quad {\rm with}\quad
\theta_1,\dots,\theta_m\in\mathbb R.
\end{displaymath}
Then,
\begin{eqnarray*}
 \nabla\gamma_N(J) & = &
 \left(\frac{1}{NT}\sum_{i = 1}^{N}\left(
 2\sum_{\ell = 1}^{m}\theta_{\ell}\int_{0}^{T}
 \varphi_j(X_{s}^{i})\varphi_{\ell}(X_{s}^{i})ds
 - 2\int_{0}^{T}\varphi_j(X_{s}^{i})dX_{s}^{i}\right)\right)_{j\in\{1,\dots,m\}}\\
 & = &
 2(\widehat{\bf\Psi}_m(\theta_1,\dots,\theta_m)^* -\widehat{\bf X}_m)
\end{eqnarray*}
where
\begin{displaymath}
\widehat{\bf\Psi}_m :=
\left(\frac{1}{NT}\sum_{i = 1}^{N}\int_{0}^{T}
\varphi_j(X_{s}^{i})\varphi_{\ell}(X_{s}^{i})ds\right)_{j,\ell\in\{1,\dots,m\}}
\end{displaymath}
and
\begin{displaymath}
\widehat{\bf X}_m :=\left(\frac{1}{NT}
\sum_{i = 1}^{N}\int_{0}^{T}
\varphi_j(X_{s}^{i})dX_{s}^{i}\right)_{j\in\{1,\dots,m\}}.
\end{displaymath}
The symmetric matrix $\widehat{\bf\Psi}_m$ is positive semidefinite because
\begin{displaymath}
\mathbf u^*\widehat{\bf\Psi}_m\mathbf u =
\frac{1}{NT}\sum_{i = 1}^{N}\int_{0}^{T}
\left(\sum_{j = 1}^{m}\mathbf u_j\varphi_j(X_{s}^{i})\right)^2ds\geqslant 0
\end{displaymath}
for every $\mathbf u\in\mathbb R^m$. If in addition $\widehat{\bf\Psi}_m$ is invertible, it is positive definite, and then
\begin{equation}\label{projection_LS_estimator_expression}
\widehat b_m =
\sum_{j = 1}^{m}\widehat\theta_j\varphi_j
\quad {\rm with}\quad
\widehat{\bf\theta}_m :=
(\widehat\theta_1,\dots,\widehat\theta_m)^* =
\widehat{\bf\Psi}_{m}^{-1}\widehat{\bf X}_m
\end{equation}
is the only minimizer of $\gamma_N(.)$ on $\mathcal S_m$, called the projection least squares estimator of $b$.
\\
\\
{\bf Remarks:}
\begin{enumerate}
 \item $\widehat{\bf\Psi}_m = (\langle\varphi_j,\varphi_{\ell}\rangle_N)_{j,\ell}$, where
 \begin{displaymath}
 \langle\varphi,\psi\rangle_N :=
 \frac{1}{NT}\sum_{i = 1}^{N}\int_{0}^{T}
 \varphi(X_{s}^{i})\psi(X_{s}^{i})ds
 \end{displaymath}
 for every continuous functions $\varphi,\psi :\mathbb R\rightarrow\mathbb R$.
 \item $\widehat{\bf X}_m = (\langle b,\varphi_j\rangle_N)_{j}^{*} +\widehat{\bf E}_m$, where
 \begin{displaymath}
 \widehat{\bf E}_m :=
 \left(\frac{1}{NT}\sum_{i = 1}^{N}\int_{0}^{T}
 \varphi_j(X_{s}^{i})(
 \sigma(X_{s}^{i})dB_{s}^{i} +
 \gamma(X_{s}^{i})d\mathfrak Z_{s}^{i})\right)_{j\in\{1,\dots,m\}}^{*}.
 \end{displaymath}
\end{enumerate}
Let us introduce the two following deterministic matrices related to the previous random ones:
\begin{displaymath}
{\bf\Psi}_{m,\sigma} := NT\mathbb E(\widehat{\bf E}_m\widehat{\bf E}_{m}^{*})
\quad {\rm and}\quad
{\bf\Psi}_m :=
\mathbb E(\widehat{\bf\Psi}_m) =
(\langle\varphi_j,\varphi_{\ell}\rangle_{f_T})_{j,\ell},
\end{displaymath}
where $\langle .,.\rangle_{f_T}$ is the usual scalar product in $\mathbb L^2(I,f_T(x)dx)$. The following lemma will be crucial in Section \ref{section_risk_bound} to evaluate the order of the variance term in the risk bound on our projection least squares estimator of $b$.
%

% Lemma : Preliminary bound on the variance term.

%
\begin{lemma}\label{preliminary_lemma_variance}
${\rm trace}(\mathbf\Psi_{m}^{-1/2}\mathbf\Psi_{m,\sigma}\mathbf\Psi_{m}^{-1/2})\leqslant (\|\sigma\|_{\infty}^{2} +\lambda\mathfrak c_{\zeta^2}\|\gamma\|_{\infty}^{2})m$.
\end{lemma}
Let us conclude this section with few words about the extension of the projection least squares estimation for multidimensional diffusion processes, and on the reason why the probability distribution of $X_s$, $s\in (0,T]$, needs to have a density with respect to Lebesgue's measure with a sharp bound as (\ref{sufficient_condition_density_1}).
\\
\\
{\bf Remarks:}
\begin{enumerate}
 \item Assume that $X$ is a $d$-dimensional diffusion process with $d\in\mathbb N^*$. A natural extension of the objective function $\gamma_N$ would be defined by
 \begin{displaymath}
 \gamma_{d,N}(\tau) :=
 \frac{1}{NT}\sum_{i = 1}^{N}\left(
 \int_{0}^{T}\|\tau(X_{s}^{i})\|_{2,d}^{2}ds - 2\int_{0}^{T}\langle\tau(X_{s}^{i}),dX_{s}^{i}\rangle_{2,d}\right)
 \textrm{ $;$ }
 \tau\in\mathcal S_{d,m},
 \end{displaymath}
 where $\mathcal S_{d,m} := {\rm span}(\{\varphi_{j_1}\otimes\dots\otimes\varphi_{j_d}\textrm{ $;$ }j_1,\dots,j_d\in\{1,\dots,m\}\})^d$. This is out of the scope of our paper, but note that one could get an expression of the solution of the minimization problem $\min_{\mathcal S_{d,m}}\gamma_{d,N}$ similar to (\ref{projection_LS_estimator_expression}) but involving hypermatrices in the spirit of Dussap (2021). So, except in the very special case where the components of $X$ are independent, to extend our estimation method to the multidimensional framework is not straightforward.
 \item By the Fubini-Tonelli theorem, for every continuous functions $\varphi,\psi :\mathbb R\rightarrow\mathbb R$,
 \begin{displaymath}
 \mathbb E(\langle\varphi,\psi\rangle_N) =
 \int_{-\infty}^{\infty}\varphi(x)\psi(x)\mathbb P_T(dx)
 \quad {\rm with}\quad
 \mathbb P_T(dx) :=\frac{1}{T}\int_{0}^{T}\mathbb P_{X_s}(dx)ds,
 \end{displaymath}
 and then one could think that there is no need to assume that $\mathbb P_{X_s}(dx) = p_s(x_0,x)dx$ for any $s\in (0,T]$. However, since the drift function $b$ may be unbounded, a sharp bound on $\mathbb P_{X_s}(dx)$ as (\ref{sufficient_condition_density_1}) is required to show that
 \begin{displaymath}
 \int_{-\infty}^{\infty}b(x)^4\mathbb P_T(dx) <\infty
 \end{displaymath}
 as assumed in the beginning of Section \ref{section_objective_function}. Note also that if the Malliavin covariance (matrix) of $X_s$, $s\in (0,T]$, is almost surely invertible (Bouleau-Hirsch's condition), then the probability distribution of $X_s$ has a density with respect to Lebesgue's measure, but not necessarily with a sharp bound as (\ref{sufficient_condition_density_1}). For instance, $X_s$ satisfies the Bouleau-Hirsch condition under Assumption 2.4 on $b$, $\sigma$ and $\gamma$ in Bichteler and Jacod (1983).
\end{enumerate}
%

% Section : Risk bound on the projection least squares estimator.

%
\section{Risk bound on the projection least squares estimator}\label{section_risk_bound}
This section deals with a risk bound on the truncated estimator
\begin{displaymath}
\widetilde b_m :=\widehat b_m\mathbf 1_{\Lambda_m},
\end{displaymath}
where
\begin{displaymath}
\Lambda_m :=
\left\{L(m)(\|\widehat{\bf\Psi}_{m}^{-1}\|_{\rm op}\vee 1)
\leqslant\mathfrak c_T\frac{NT}{\log(NT)}
\right\}
\end{displaymath}
with
\begin{displaymath}
L(m) =
1\vee\left(
\sup_{x\in I}\sum_{j = 1}^{m}\varphi_j(x)^2\right)
\quad {\rm and}\quad
\mathfrak c_T =
\frac{3\log(3/2) - 1}{8T}.
\end{displaymath}
On the event $\Lambda_m$, $\widehat{\bf\Psi}_m$ is invertible because
\begin{displaymath}
\inf\{{\rm sp}(\widehat{\bf\Psi}_m)\}
\geqslant\frac{L(m)}{\mathfrak c_T}\cdot\frac{\log(NT)}{NT},
\end{displaymath}
and then $\widetilde b_m$ is well-defined. In the sequel, $m$ fulfills the following assumption.
%

% Assumption : Condition on m.

%
\begin{assumption}\label{condition_m}
$\displaystyle{L(m)(\|\mathbf\Psi_{m}^{-1}\|_{\rm op}\vee 1)
\leqslant\frac{\mathfrak c_T}{2}\cdot\frac{NT}{\log(NT)}}$.
\end{assumption}
The above condition is a generalization of the so-called {\it stability condition} introduced in the nonparametric regression framework in Cohen et al. (2013), and already extended to the independent copies of continuous diffusion processes framework in Comte and Genon-Catalot (2020b).
%

% Theorem : Risk bound.

%
\begin{theorem}\label{risk_bound}
Under Assumption \ref{condition_m}, there exists a constant $\mathfrak c_{\ref{risk_bound}} > 0$, not depending on $m$ and $N$, such that
\begin{equation}\label{risk_bound_1}
\mathbb E(\|\widetilde b_m - b_I\|_{N}^{2})
\leqslant
\min_{\tau\in\mathcal S_m}\|\tau - b_I\|_{f_T}^{2} +
\mathfrak c_{\ref{risk_bound}}\left(\frac{m}{NT} +
\frac{1}{N}\right).
\end{equation}
\end{theorem}
\noindent
{\bf Remarks:}
\begin{itemize}
 \item Note that Theorem \ref{risk_bound} says first that the bound on the variance of $\widetilde b_m$ is of order $m/N$ as in the nonparametric regression framework, in which it is optimal (see Comte and Genon-Catalot (2020a), Theorem 1). For this reason, (\ref{risk_bound_1}) should be near to optimality, but it's a difficult challenge, out of the scope of this paper, to establish a lower bound on $\widehat b_m$. Note that even for continuous diffusion processes, no lower bound has been established (see Comte and Genon-Catalot (2020b)).
 \item The order of the bias term
 \begin{displaymath}
 \min_{\tau\in\mathcal S_m}
 \|\tau - b_I\|_{f_T}^{2},
 \end{displaymath}
 as well as $L(m)$ and $\|\mathbf\Psi_{m}^{-1}\|_{\rm op}$, depend on the $\varphi_j$'s. Let us evaluate them for the trigonometric basis, which is compactly supported, and for the Hermite basis, which is $\mathbb R$-supported:
 \begin{enumerate}
  \item Assume that $I = [\ell,\texttt r]$ with $\ell,\texttt r\in\mathbb R$ satisfying $\ell <\texttt r$, and that
  \begin{eqnarray*}
   \varphi_1(x) & := &
   \sqrt{\frac{1}{\texttt r -\ell}}\mathbf 1_{[\ell,\texttt r]}(x),\\
   \varphi_{2j + 1}(x) & := &
   \sqrt{\frac{2}{\texttt r -\ell}}
   \sin\left(2\pi j\frac{x -\ell}{\texttt r -\ell}\right)\mathbf 1_{[\ell,\texttt r]}(x)
   \quad {\rm and}\\
   \varphi_{2j}(x) & := &
   \sqrt{\frac{2}{\texttt r -\ell}}
   \cos\left(2\pi j\frac{x -\ell}{\texttt r -\ell}\right)\mathbf 1_{[\ell,\texttt r]}(x)
  \end{eqnarray*}
  for every $x\in [\ell,\texttt r]$ and $j\in\mathbb N^*$. On the one hand, since $\cos(.)^2 +\sin(.)^2 = 1$, there exists a constant $\mathfrak c_{\varphi} > 0$, not depending on $m$, $N$ and $T$, such that
  \begin{displaymath}
  L(m) =
  1\vee\left[\sup_{x\in I}\sum_{j = 1}^{m}\varphi_j(x)^2\right]
  \leqslant\mathfrak c_{\varphi}m.
  \end{displaymath}
  Moreover, since under the same conditions than in Remark \ref{sufficient_condition_density}, there exists $\underline{\mathfrak m} > 0$ such that $f_T(.)\geqslant\underline{\mathfrak m}$ on $I$ by Chen et al. (2017), Theorem 1.3,
  \begin{eqnarray*}
   \|\mathbf\Psi_{m}^{-1}\|_{\rm op} =
   \frac{1}{\lambda_{\min}(\mathbf\Psi_m)} & = &
   \left(\inf_{\theta :\|\theta\|_{2,m} = 1}
   \sum_{j,\ell = 1}^{m}\theta_j\theta_{\ell}[\mathbf\Psi_m]_{j,\ell}\right)^{-1}\\
   & = &
   \left[
   \inf_{\theta :\|\theta\|_{2,m} = 1}\int_{\ell}^{\tt r}\left(
   \sum_{j = 1}^{m}\theta_j\varphi_j(x)\right)^2f_T(x)dx\right]^{-1}
   \leqslant\frac{1}{\underline{\mathfrak m}}.
  \end{eqnarray*}
  Then,
  \begin{displaymath}
  m\leqslant
  \frac{1}{\mathfrak c_{\varphi}(\overline{\mathfrak m}^{-1}\vee 1)}\cdot
  \frac{\mathfrak c_T}{2}\cdot\frac{NT}{\log(NT)}
  \end{displaymath}
  satisfies Assumption \ref{condition_m}. On the other hand, consider the Fourier-Sobolev space
  \begin{displaymath}
  \mathbb W_{2}^{\beta}([\ell,\texttt r]) :=
  \left\{\varphi\in C^{\beta}([\ell,\texttt r];\mathbb R) :
  \int_{\ell}^{\tt r}\varphi^{(\beta)}(x)^2dx <\infty\right\}
  \end{displaymath}
  with $\beta\in\mathbb N^*$, and assume that $b_{|I}\in\mathbb W_{0}^{\beta}([\ell,\texttt r])$. By DeVore and Lorentz (1993), Corollary 2.4 p. 205, there exists a constant $\mathfrak c_{\beta,\ell,\texttt r} > 0$, not depending on $m$ and $N$, such that
  \begin{displaymath}
  \|\Pi_m(b_I) - b_I\|^2\leqslant
  \mathfrak c_{\beta,\ell,\texttt r}m^{-2\beta},
  \end{displaymath}
  where $\Pi_m$ is the orthogonal projection from $\mathbb L^2(I,dx)$ onto $\mathcal S_m$. Since under appropriate conditions on $b$ and $\pi$, as explained in Remark \ref{sufficient_condition_density}, $f_T$ is upper bounded on $I$ by a constant $\overline{\mathfrak m} > 0$,
  \begin{eqnarray*}
   \min_{\tau\in\mathcal S_m}
   \|\tau - b_I\|_{f_T}^{2}
   & \leqslant &
   \overline{\mathfrak m}\|\Pi_m(b_I) - b_I\|^2\\
   & \leqslant &
   \overline{\mathfrak c}_{\beta,\ell,\texttt r}m^{-2\beta}
   \quad {\rm with}\quad
   \overline{\mathfrak c}_{\beta,\ell,\texttt r} =
   \mathfrak c_{\beta,\ell,\texttt r}
   \overline{\mathfrak m}.
  \end{eqnarray*}
  In conclusion, by Theorem \ref{risk_bound}, there exists a constant $\overline{\mathfrak c}_{\ref{risk_bound},1} > 0$, not depending $m$ and $N$, such that
  \begin{displaymath}
  \mathbb E(\|\widetilde b_m - b_I\|_{N}^{2})\leqslant
  \overline{\mathfrak c}_{\ref{risk_bound},1}\left(
  m^{-2\beta} +\frac{m}{N}\right),
  \end{displaymath}
  and then the bias-variance tradeoff is reached by (the risk bound on) $\widetilde b_m$ for $m$ of order $N^{1/(1 + 2\beta)}$.
  \item Assume that $\varphi_j = h_{j - 1}$ for any $j\in\{1,\dots,m\}$, where $(h_n)_{n\in\mathbb N}$ is the Hermite basis: for every $x\in\mathbb R$ and $n\in\mathbb N$,
  \begin{displaymath}
  h_n(x) :=
  (2^nn!\sqrt\pi)^{-1/2}H_n(x)e^{-x^2/2}
  \quad {\rm with}\quad
  H_n(x) = (-1)^ne^{x^2}\frac{d^n}{dx^n}e^{-x^2}.
  \end{displaymath}
  On the one hand, by Abramowitz and Stegun (1964), $\|\varphi_j\|_{\infty}\leqslant\pi^{-1/4}$ and then $L(m)\leqslant m$. Moreover, assume that the conditions (1) and (2) in Remark \ref{sufficient_condition_density} are satisfied by $b$ and $\pi$ with $\alpha\in [1,2)$. By Chen et al. (2017), Theorem 1.3, for every $s\in (0,T]$, there exist two constants $\overline{\mathfrak c}_p,\overline{\mathfrak m}_p > 0$, not depending on $s$, such that for any $x\in\mathbb R$,
  \begin{eqnarray*}
   p_s(x_0,x) & \geqslant & \overline{\mathfrak c}_p\left[
   s^{-1/2}\exp\left(-\overline{\mathfrak m}_p\frac{(x - x_0)^2}{s}\right) +\frac{s}{(s^{1/2} +
   |x - x_0|)^{1 +\alpha}}\right]\\
   & \geqslant &
   \frac{\overline{\mathfrak c}_ps}{(2T + 4x_{0}^{2} + 4x^2)^{(1 +\alpha)/2}}
   \geqslant
   \frac{\overline{\mathfrak c}_ps}{[4\vee (2T + 4x_{0}^{2})]^{(1 +\alpha)/2}}\cdot
   \frac{1}{(1 + x^2)^{(1 +\alpha)/2}}.
  \end{eqnarray*}
  So,
  \begin{displaymath}
  f_T(x)\geqslant
  \frac{\mathfrak c_{f_T}}{(1 + x^2)^{(1 +\alpha)/2}}
  \quad {\rm with}\quad
  \mathfrak c_{f_T} =
  \frac{\overline{\mathfrak c}_pT}{2[4\vee (2T + 4x_{0}^{2})]^{(1 +\alpha)/2}}.
  \end{displaymath}
  Since $\alpha\in [1,2)$, by Comte and Genon-Catalot (2020a), Proposition 9, there exists a constant $\mathfrak c_{\bf\Psi} > 0$, not depending on $m$, such that
  \begin{displaymath}
  \|\mathbf\Psi_{m}^{-1}\|_{\rm op}\leqslant\mathfrak c_{\bf\Psi}m^{(1 +\alpha)/2}.
  \end{displaymath}
  Then,
  \begin{displaymath}
  m\leqslant
  \left[\frac{1}{\mathfrak c_{\bf\Psi}\vee 1}\cdot
  \frac{\mathfrak c_T}{2}\cdot\frac{NT}{\log(NT)}\right]^{2/(3 +\alpha)}
  \end{displaymath}
  satisfies Assumption \ref{condition_m}. On the other hand, assume that $b\in\mathbb W_{\tt H}^{\beta}(D)$, where $\mathbb W_{\tt H}^{\beta}(D)$ is the Hermite-Sobolev ball defined by
  \begin{displaymath}
  \mathbb W_{\tt H}^{\beta}(D) :=
  \left\{\varphi\in\mathbb L^2(\mathbb R) :
  \sum_{n = 0}^{\infty}
  n^{\beta}\langle\varphi,h_n\rangle^2\leqslant D\right\}
  \end{displaymath}
  with $\beta > (3 +\alpha)/2 - 1$ and $D > 0$. Then, $\|\Pi_m(b) - b\|^2\leqslant Dm^{-\beta}$ (see Belomestny et al. (2019), Section 4.2), and since $f_T$ is upper bounded on $\mathbb R$ by a constant $\overline{\mathfrak m} > 0$ (see Remark \ref{sufficient_condition_density}),
  \begin{eqnarray*}
   \min_{\tau\in\mathcal S_m}
   \|\tau - b\|_{f_T}^{2}
   & \leqslant &
   \overline{\mathfrak m}\|\Pi_m(b) - b\|^2\\
   & \leqslant &
   \mathfrak c_{\tt H}m^{-\beta}
   \quad {\rm with}\quad
   \mathfrak c_{\tt H} =
   D\overline{\mathfrak m}.
  \end{eqnarray*}
  In conclusion, by Theorem \ref{risk_bound}, there exists a constant $\overline{\mathfrak c}_{\ref{risk_bound},2} > 0$, not depending $m$ and $N$, such that
  \begin{displaymath}
  \mathbb E(\|\widetilde b_m - b\|_{N}^{2})\leqslant
  \overline{\mathfrak c}_{\ref{risk_bound},2}\left(
  m^{-\beta} +\frac{m}{N}\right),
  \end{displaymath}
  and then the bias-variance tradeoff is reached by (the risk bound on) $\widetilde b_m$ for $m$ of order $N^{1/(1 +\beta)}$.
 \end{enumerate}
\end{itemize}
%

% Section : Model selection.

%
\section{Model selection}\label{section_model_selection}
First of all, let us state some additional assumptions on the $\varphi_j$'s, on the L\'evy measure $\pi_{\lambda}$ and on the density function $f_T$.
%

% Assumption : Nested models.

%
\begin{assumption}\label{nested_models}
The spaces $\mathcal S_m$, $m\in\{1,\dots,N_T\}$, satisfy the following conditions:
\begin{enumerate}
 \item There exists a constant $\mathfrak c_{\varphi}\geqslant 1$, not depending on $N$, such that for every $m\in\{1,\dots,N_T\}$, the basis $(\varphi_1,\dots,\varphi_m)$ of $\mathcal S_m$ satisfies
 \begin{displaymath}
 L(m) = 1\vee\left[\sup_{x\in I}
 \sum_{j = 1}^{m}\varphi_j(x)^2\right]
 \leqslant\mathfrak c_{\varphi}^2m.
 \end{displaymath}
 \item (Nested spaces) For every $m,m'\in\{1,\dots,N_T\}$, if $m>m'$, then $\mathcal S_{m'}\subset\mathcal S_m$.
\end{enumerate}
\end{assumption}
The proof of the oracle inequality satisfied by the adaptive estimator relies on a Bernstein type inequality (see Lemma \ref{Bernstein_type_inequality}). In order to control the {\it big jumps} of $X$ in the proof of Lemma \ref{Bernstein_type_inequality}, the L\'evy measure $\pi_{\lambda}$ needs to fulfill the following assumption.
%

% Assumption : Sub-exponential Lévy measure.

%
\begin{assumption}\label{sub_exponential_Levy_measure}
The L\'evy measure $\pi_{\lambda}$ is sub-exponential: there exist positive constants $\mathfrak a,\mathfrak b > 0$ such that, for every $x > 1$,
\begin{displaymath}
\pi_{\lambda}((-x,x)^c)
\leqslant\mathfrak ae^{-\mathfrak b|x|}.
\end{displaymath}
\end{assumption}
%

% Assumption : Bounded density.

%
\begin{assumption}\label{bounded_density}
The density function $f_T$ is bounded.
\end{assumption}
Note that Remark \ref{sufficient_condition_density} provides a sufficient condition for Assumption \ref{bounded_density} to be satisfied by $f_T$. Now, let us introduce the set
\begin{displaymath}
\widehat{\mathcal M}_N :=\left\{m\in\{1,\dots,N_T\} :
\mathfrak c_{\varphi}^{2}m(\|\widehat{\bf\Psi}_{m}^{-1}\|_{\rm op}^{2}\vee 1)
\leqslant\mathfrak d_T\frac{NT}{\log(NT)}\right\}
\end{displaymath}
with
\begin{displaymath}
\mathfrak d_T =
\min\left\{\frac{\mathfrak c_T}{2},
\frac{1}{64\mathfrak c_{\varphi}^{2}T(\|f_T\|_{\infty} +\sqrt{\mathfrak c_T/2}/(3\mathfrak c_{\varphi}))}
\right\},
\end{displaymath}
as well as its theoretical counterpart
\begin{displaymath}
\mathcal M_N :=\left\{m\in\{1,\dots,N_T\} :
\mathfrak c_{\varphi}^{2}m(\|\mathbf\Psi_{m}^{-1}\|_{\rm op}^{2}\vee 1)
\leqslant\frac{\mathfrak d_T}{4}\cdot\frac{NT}{\log(NT)}\right\}.
\end{displaymath}
Let us also consider
\begin{displaymath}
\widehat m :=\arg\min_{m\in\widehat{\mathcal M}_N}
\{-\|\widehat b_m\|_{N}^{2} + {\rm pen}(m)\},
\end{displaymath}
where
\begin{displaymath}
{\rm pen}(m) :=\mathfrak c_{\rm cal}\frac{m}{NT}
\textrm{ $;$ }\forall m\in\{1,\dots,N_T\}
\end{displaymath}
and $\mathfrak c_{\rm cal} > 0$ is a deterministic constant to calibrate.
%

% Theorem : Model selection theorem.

%
\begin{theorem}\label{th_model_selection}
Under Assumptions \ref{condition_m}, \ref{nested_models}, \ref{sub_exponential_Levy_measure} and \ref{bounded_density}, there exists a constant $\mathfrak c_{\ref{th_model_selection}} > 0$, not depending on $N$, such that
\begin{displaymath}
\mathbb E(\|\widehat b_{\widehat m} - b_I\|_{N}^{2})
\leqslant\mathfrak c_{\ref{th_model_selection}}\left[
\min_{m\in\mathcal M_N}\left\{
\min_{\tau\in\mathcal S_m}\|\tau - b_I\|_{f_T}^{2} +\frac{m}{NT}\right\} +
\frac{1}{N}\right].
\end{displaymath}
\end{theorem}
Note that Theorem \ref{th_model_selection} provides a risk bound on the adaptive estimator $\widehat b_{\widehat m}$ of same order than the minimal risk bound on $\widetilde b_m$ (see Theorem \ref{risk_bound}) for $m$ taken in $\mathcal M_N$.
%

% Section : Numerical experiments.

%
\section{Numerical experiments}\label{section_numerical_experiments}
First of all, let us recall that all figures and tables related to this section are postponed to Appendix \ref{section_tables_figures}.

\begin{figure}[!h]
\centering
\includegraphics[scale=0.4]{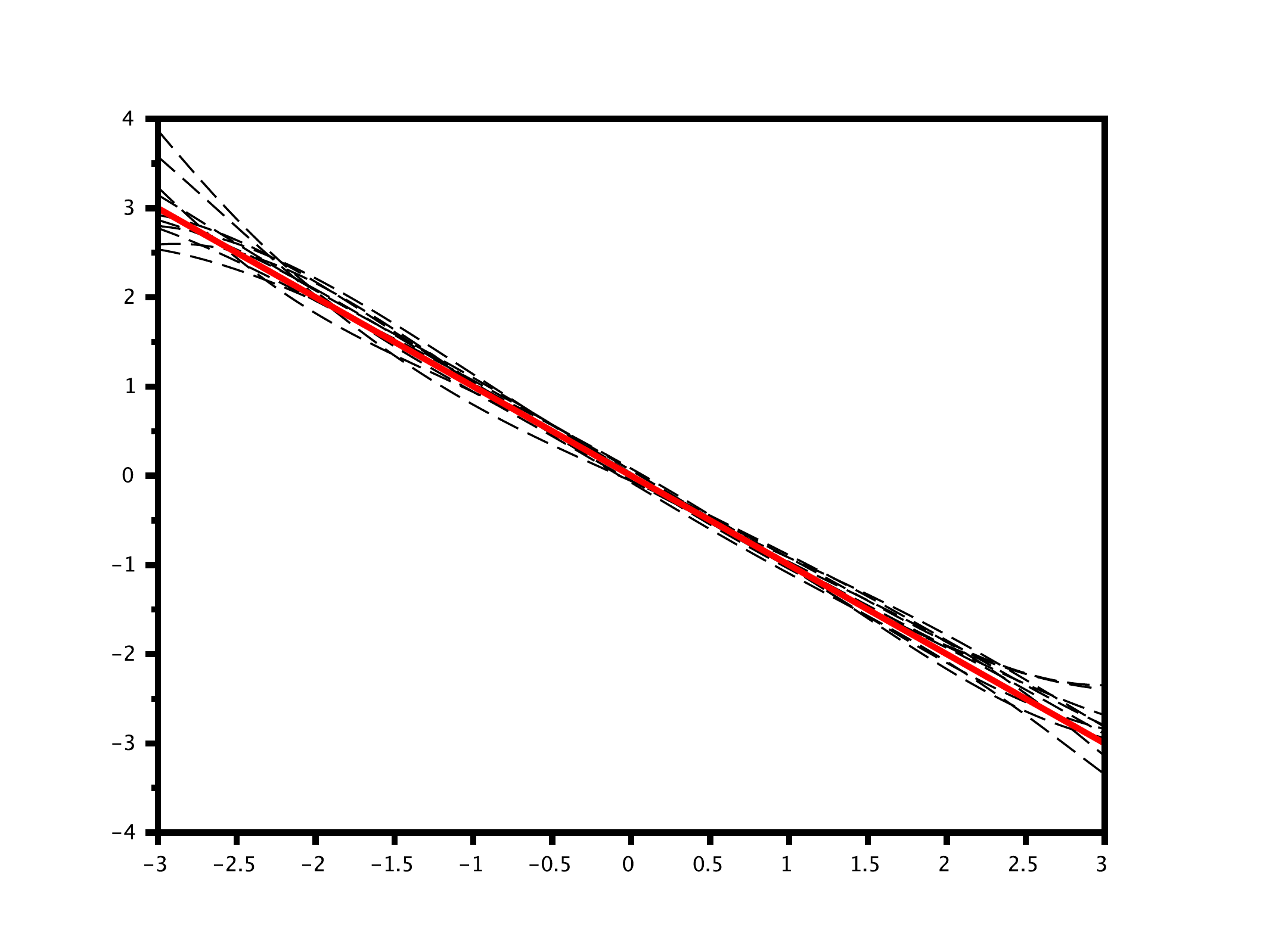} 
\caption{Plots of $b$ and of $10$ adaptive estimations for Model 1 ($\overline{\widehat m} = 5.3$).}
\label{Model_1_plot}
\end{figure}

\begin{figure}[!h]
\centering
\includegraphics[scale=0.4]{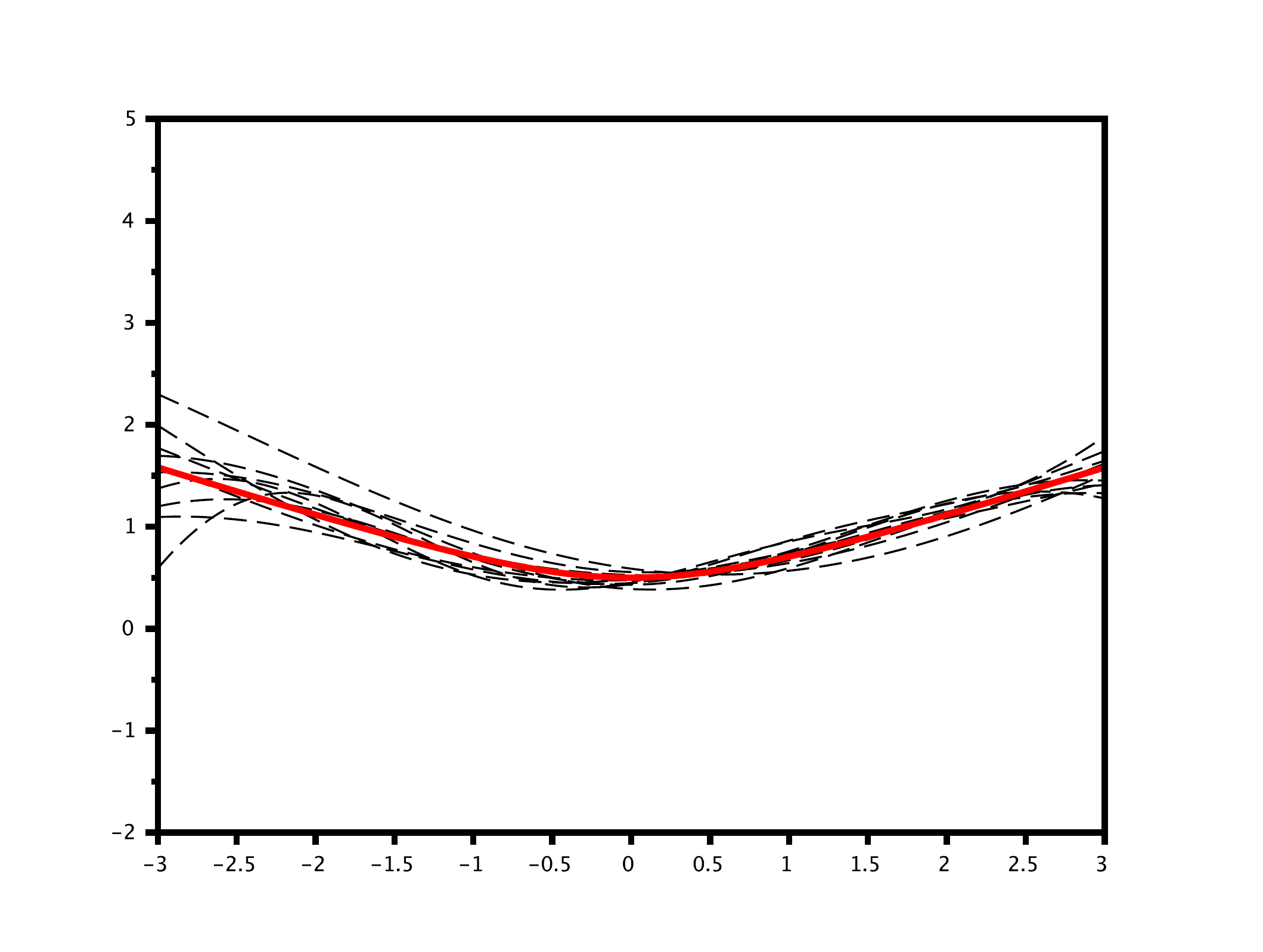} 
\caption{Plots of $b$ and of $10$ adaptive estimations for Model 2 ($\overline{\widehat m} = 4.2$).}
\label{Model_2_plot}
\end{figure}

\begin{figure}[!h]
\centering
\includegraphics[scale=0.4]{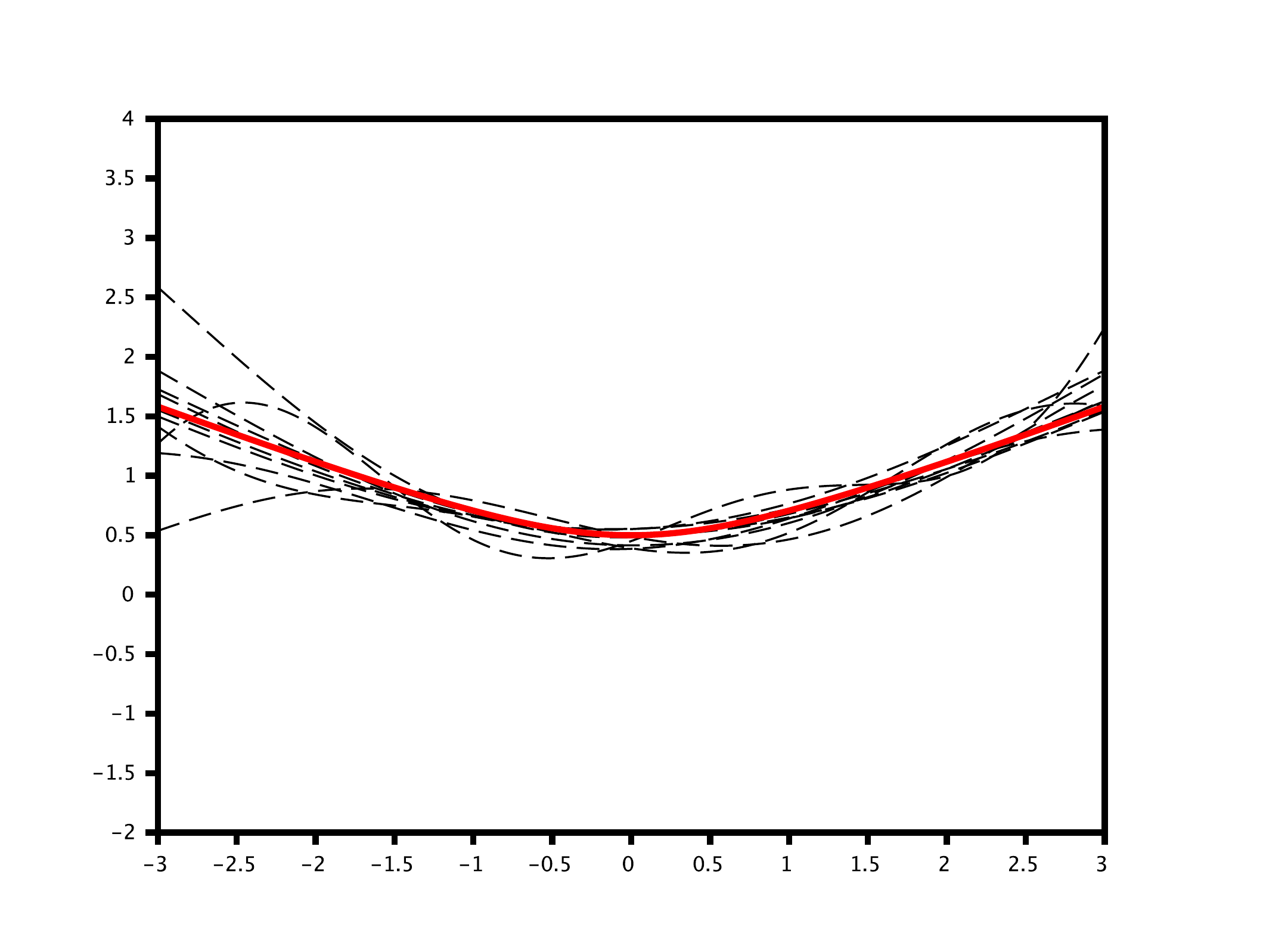} 
\caption{Plots of $b$ and of $10$ adaptive estimations for Model 3 ($\overline{\widehat m} = 4.1$).}
\label{Model_3_plot}
\end{figure}

Some numerical experiments on our estimation method of $b$ in Equation (\ref{main_equation}) are presented in this subsection when the common distribution $\pi$ of the $\zeta_n$'s is standard normal, which implies that $\mathfrak Z = Z$, and the intensity of the (usual) Poisson process $\nu$ is $\lambda = 0.5$. The estimation method investigated on the theoretical side in Sections \ref{section_risk_bound} and \ref{section_model_selection} is implemented here for the three following models:
\begin{enumerate}
 \item $\displaystyle{X_t = 0.5 -\int_{0}^{t}X_sds + 0.5B_t + Z_t}$ (linear/additive noise).
 \item $\displaystyle{X_t = 0.5 + 0.5\int_{0}^{t}\sqrt{1 + X_{s}^{2}}ds + 0.5B_t + Z_t}$ (nonlinear/additive noise).
 \item $\displaystyle{X_t = 0.5 + 0.5\int_{0}^{t}\sqrt{1 + X_{s}^{2}}ds + 0.5\int_{0}^{t}(1 +\cos(X_s)^2)dB_s + Z_t}$ (nonlinear/multiplicative noise).
\end{enumerate}
For each of the three previous models, our adaptive estimator of $b$ is computed on $I = [-3,3]$ from $N = 400$ paths of the process $X$ observed along the dissection $\{lT/n\textrm{ $;$ }l = 1,\dots,n\}$ of $[0,T]$ with $n = 200$ and $T = 5$, when $(\varphi_1,\dots,\varphi_m)$ is the $m$-dimensional trigonometric basis for every $m\in\{1,\dots,6\}$. This experiment is repeated $100$ times, and the means and the standard deviations of the MISE of $\widehat b_{\widehat m}$ are stored in Table \ref{table_MISE}. Moreover, for each model, $10$ estimations (dashed black curves) of $b$ (red curve) are plotted on Figures \ref{Model_1_plot}, \ref{Model_2_plot} and \ref{Model_3_plot}.

On the one hand, on average, the MISE of our adaptive estimator is slightly increasing with the {\it complexity} of the model: 0.1251 (Model 1) < 0.1469 (Model 2) < 0.1825 (Model 3). The same remark holds for its standard deviation: 0.0950 (Model 1) < 0.1688 (Model 2) < 0.1928 (Model 3). This means that the more the model is {\it complex}, the more the quality of the estimation degrades, and it's visible on Figures \ref{Model_1_plot}, \ref{Model_2_plot} and \ref{Model_3_plot}. However, for all the three previous models, the MISE of $\widehat b_{\widehat m}$ remains good and of same order than for models without jump component (see Comte and Genon-Catalot (2020b), Section 4). On the other hand, the model selection procedure works well because it doesn't systematically select the lowest or the highest possible value of $m$ (see $\overline{\widehat m}$ on Figures \ref{Model_1_plot}, \ref{Model_2_plot} and \ref{Model_3_plot}). For Model 1, the procedure is stable because $\widehat m$ has a low estimated standard deviation (0.4830). The procedure remains satisfactory for Models 2 and 3, with ${\rm StD}(\widehat m) = 1.1353$ and ${\rm StD}(\widehat m) = 1.2867$ respectively, but not as much as for Model 1.
\appendix
%

% Section : Proofs

%
\section{Proofs}\label{section_proofs}
%

% Subsection : Proof of Lemma preliminary_lemma_variance.

%
\subsection{Proof of Lemma \ref{preliminary_lemma_variance}}
First, let us show that the symmetric matrix $\mathbf\Psi_{m,\sigma}$ is positive semidefinite. Indeed, for any $y\in\mathbb R^m$,
\begin{eqnarray*}
 y^*\mathbf\Psi_{m,\sigma}y & = &
 \frac{1}{NT}\sum_{j,\ell = 1}^{m}y_jy_{\ell}
 \sum_{i,k = 1}^{N}\mathbb E\left(
 \left(\int_{0}^{T}\varphi_j(X_{s}^{i})(\sigma(X_{s}^{i})dB_{s}^{i} +
 \gamma(X_{s}^{i})d\mathfrak Z_{s}^{i})\right)\right.\\
 & &
 \hspace{5cm}\left.\times
 \left(\int_{0}^{T}\varphi_{\ell}(X_{s}^{k})(\sigma(X_{s}^{k})dB_{s}^{k} +
 \gamma(X_{s}^{k})d\mathfrak Z_{s}^{k})\right)\right)\\
 & = &
 \frac{1}{NT}\mathbb E\left[\left(
 \sum_{i = 1}^{N}\int_{0}^{T}\tau_y(X_{s}^{i})(\sigma(X_{s}^{i})dB_{s}^{i} +
 \gamma(X_{s}^{i})d\mathfrak Z_{s}^{i})\right)^2\right]
 \geqslant 0
 \quad {\rm with}\quad
 \tau_y(.) :=\sum_{j = 1}^{m}y_j\varphi_j(.).
\end{eqnarray*}
Moreover, by the isometry property of It\^o's integral (with respect to $B$), by the isometry type property of the stochastic integral with respect to $\mathfrak Z$, and since $\sigma$ and $\gamma$ are bounded, for every $j\in\{1,\dots,m\}$,
\begin{eqnarray}
 y^*\mathbf\Psi_{m,\sigma}y
 & = &
 \frac{1}{T}\mathbb E\left[\left(\int_{0}^{T}\tau_y(X_s)\sigma(X_s)dB_s\right)^2\right] +
 \frac{1}{T}\mathbb E\left[
 \left(\int_{0}^{T}\tau_y(X_s)\gamma(X_s)d\mathfrak Z_s\right)^2\right]
 \nonumber\\
 & = &
 \frac{1}{T}\int_{0}^{T}\mathbb E(\tau_y(X_s)^2\sigma(X_s)^2)ds
 +\frac{\lambda\mathfrak c_{\zeta^2}}{T}
 \int_{0}^{T}\mathbb E(\tau_y(X_s)^2\gamma(X_s)^2)ds
 \nonumber\\
 \label{preliminary_lemma_variance_1}
 & \leqslant &
 (\|\sigma\|_{\infty}^{2} +\lambda\mathfrak c_{\zeta^2}\|\gamma\|_{\infty}^{2})
 \int_{-\infty}^{\infty}
 \left(\sum_{j = 1}^{m}y_j\varphi_j(x)\right)^2f_T(x)dx =
 (\|\sigma\|_{\infty}^{2} +\lambda\mathfrak c_{\zeta^2}\|\gamma\|_{\infty}^{2})
 \|\mathbf\Psi_{m}^{1/2}y\|_{2,m}^{2}.
\end{eqnarray}
Therefore, since $\mathbf\Psi_{m,\sigma}$ is positive semidefinite, and by Inequality (\ref{preliminary_lemma_variance_1}),
\begin{eqnarray*}
 {\rm trace}(\mathbf\Psi_{m}^{-1/2}\mathbf\Psi_{m,\sigma}\mathbf\Psi_{m}^{-1/2})
 & \leqslant &
 m\|\mathbf\Psi_{m}^{-1/2}\mathbf\Psi_{m,\sigma}\mathbf\Psi_{m}^{-1/2}\|_{\rm op}\\
 & = &
 m\cdot
 \sup\{y^*\mathbf\Psi_{m,\sigma}y\textrm{ $;$ }
 y\in\mathbb R^m\textrm{ and }\|\mathbf\Psi_{m}^{1/2}y\|_{2,m} = 1\}\\
 & \leqslant &
 (\|\sigma\|_{\infty}^{2} +\lambda\mathfrak c_{\zeta^2}\|\gamma\|_{\infty}^{2})m.
\end{eqnarray*}
\hfill $\Box$
%

% Subsection : Proof of Theorem risk_bound.

%
\subsection{Proof of Theorem \ref{risk_bound}}
The proof of Theorem \ref{risk_bound} relies on the two following lemmas.
%

% Lemma : Bound on the variance of the remainder term.

%
\begin{lemma}\label{bound_variance_remainder}
There exists a constant $\mathfrak c_{\ref{bound_variance_remainder}} > 0$, not depending on $m$ and $N$, such that
\begin{displaymath}
\mathbb E(|\widehat{\bf E}_{m}^{*}\widehat{\bf E}_m|^2)
\leqslant
\mathfrak c_{\ref{bound_variance_remainder}}
\frac{mL(m)^2}{N^2}.
\end{displaymath}
\end{lemma}
%

% Lemma : Deviation probabilities.

%
\begin{lemma}\label{deviation_probabilities}
Consider the event
\begin{displaymath}
\Omega_m :=
\left\{\sup_{\tau\in\mathcal S_m}
\left|\frac{\|\tau\|_{N}^{2}}{\|\tau\|_{f_T}^{2}} - 1\right|\leqslant\frac{1}{2}
\right\}.
\end{displaymath}
Under Assumptions \ref{condition_m}, there exists a constant $\mathfrak c_{\ref{deviation_probabilities}} > 0$, not depending on $m$ and $N$, such that
\begin{displaymath}
\mathbb P(\Omega_{m}^{c})
\leqslant
\frac{\mathfrak c_{\ref{deviation_probabilities}}}{N^7}
\quad\textrm{and}\quad
\mathbb P(\Lambda_{m}^{c})
\leqslant
\frac{\mathfrak c_{\ref{deviation_probabilities}}}{N^7}.
\end{displaymath}
\end{lemma}
The proof of Lemma \ref{bound_variance_remainder} is postponed to Subsubsection \ref{proof_variance_remainder}, and the proof of Lemma \ref{deviation_probabilities} remains the same as the proof of Comte and Genon-Catalot (2020b), Lemma 6.1, because $(B^1,Z^1),\dots,(B^N,Z^N)$ are independent.
%

% Subsubsection : Steps of the proof.

%
\subsubsection{Steps of the proof}
First of all,
\begin{eqnarray*}
 \|\widetilde b_m - b_I\|_{N}^{2} & = &
 \|b_I\|_{N}^{2}\mathbf 1_{\Lambda_{m}^{c}} +
 \|\widehat b_m - b_I\|_{N}^{2}\mathbf 1_{\Lambda_m}\\
 & = &
 \|b_I\|_{N}^{2}\mathbf 1_{\Lambda_{m}^{c}} +
 \|\widehat b_m - b_I\|_{N}^{2}\mathbf 1_{\Lambda_m\cap\Omega_m} +
 \|\widehat b_m - b_I\|_{N}^{2}\mathbf 1_{\Lambda_m\cap\Omega_{m}^{c}} =:
 U_1 + U_2 + U_3.
\end{eqnarray*}
Let us find suitable bounds on $\mathbb E(U_1)$, $\mathbb E(U_2)$ and $\mathbb E(U_3)$.
\begin{itemize}
 \item {\bf Bound on $\mathbb E(U_1)$.} By Cauchy-Schwarz's inequality,
 \begin{eqnarray*}
  \mathbb E(U_1) & \leqslant &
  \mathbb E(\|b_I\|_{N}^{4})^{1/2}\mathbb P(\Lambda_{m}^{c})^{1/2}
  \leqslant
  \mathbb E\left(\frac{1}{T}\int_{0}^{T}b_I(X_t)^4dt\right)^{1/2}
  \mathbb P(\Lambda_{m}^{c})^{1/2}\\
  & \leqslant &
  \mathfrak c_1\mathbb P(\Lambda_{m}^{c})^{1/2} <\infty
  \quad {\rm with}\quad
  \mathfrak c_1 =\left(\int_{-\infty}^{\infty}b_I(x)^4f_T(x)dx\right)^{1/2} <\infty.
 \end{eqnarray*}
 \item {\bf Bound on $\mathbb E(U_2)$.} Let $\Pi_{N,m}(.)$ be the orthogonal projection from $\mathbb L^2(I,f_T(x)dx)$ onto $\mathcal S_m$ with respect to the empirical scalar product $\langle .,.\rangle_N$. Then,
 \begin{equation}\label{risk_bound_2}
 \|\widehat b_m - b_I\|_{N}^{2} =
 \|\widehat b_m -\Pi_{N,m}(b_I)\|_{N}^{2} +
 \min_{\tau\in\mathcal S_m}\|b_I -\tau\|_{N}^{2}.
 \end{equation}
 As in the proof of Comte and Genon-Catalot (2020b), Proposition 2.1, on $\Omega_m$,
 \begin{displaymath}
 \|\widehat b_m -\Pi_{N,m}(b_I)\|_{N}^{2} =
 \widehat{\bf E}_{m}^{*}\widehat{\bf\Psi}_{m}^{-1}\widehat{\bf E}_m
 \leqslant 2\widehat{\bf E}_{m}^{*}\mathbf\Psi_{m}^{-1}\widehat{\bf E}_m.
 \end{displaymath}
 So,
 \begin{eqnarray*}
  \mathbb E(\|\widehat b_m -\Pi_{N,m}(b_I)\|_{N}^{2}
  \mathbf 1_{\Lambda_m\cap\Omega_m}) & \leqslant &
  2\mathbb E\left(\sum_{j,\ell = 1}^{m}[\widehat{\bf E}_m]_j[\widehat{\bf E}_m]_{\ell}
  [\mathbf\Psi_{m}^{-1}]_{j,\ell}\right)\\
  & = &
  \frac{2}{NT}\sum_{j,\ell = 1}^{m}
  [\mathbf\Psi_{m,\sigma}]_{j,\ell}
  [\mathbf\Psi_{m}^{-1}]_{j,\ell} =
  \frac{2}{NT}{\rm trace}(\mathbf\Psi_{m}^{-1/2}
  \mathbf\Psi_{m,\sigma}
  \mathbf\Psi_{m}^{-1/2}).
 \end{eqnarray*}
 Then, by Equality (\ref{risk_bound_2}) and Lemma \ref{preliminary_lemma_variance},
 \begin{eqnarray*}
  \mathbb E(U_2) & \leqslant &
  \mathbb E\left(
  \min_{\tau\in\mathcal S_m}\|b_I -\tau\|_{N}^{2}\right) +
  \frac{2}{NT}{\rm trace}(\mathbf\Psi_{m}^{-1/2}
  \mathbf\Psi_{m,\sigma}
  \mathbf\Psi_{m}^{-1/2})\\
  & \leqslant &
  \min_{\tau\in\mathcal S_m}\|b_I -\tau\|_{f_T}^{2} +
  \frac{2m}{NT}
  (\|\sigma\|_{\infty}^{2} +\lambda\mathfrak c_{\zeta^2}\|\gamma\|_{\infty}^{2}).
 \end{eqnarray*}
 \item {\bf Bound on $\mathbb E(U_3)$.} Since
 \begin{displaymath}
 \|\widehat b_m -\Pi_{N,m}(b_I)\|_{N}^{2} =
 \widehat{\bf E}_{m}^{*}\widehat{\bf\Psi}_{m}^{-1}\widehat{\bf E}_m,
 \end{displaymath}
 by the definition of the event $\Lambda_m$, and by Lemma \ref{bound_variance_remainder},
 \begin{eqnarray*}
  \mathbb E(\|\widehat b_m -\Pi_{N,m}(b_I)\|_{N}^{2}
  \mathbf 1_{\Lambda_m\cap\Omega_{m}^{c}})
  & \leqslant &
  \mathbb E(\|\widehat{\bf\Psi}_{m}^{-1}\|_{\rm op}
  |\widehat{\bf E}_{m}^{*}\widehat{\bf E}_m|
  \mathbf 1_{\Lambda_m\cap\Omega_{m}^{c}})\\
  & \leqslant &
  \frac{\mathfrak c_T}{L(m)}\cdot
  \frac{NT}{\log(NT)}
  \mathbb E(|\widehat{\bf E}_{m}^{*}\widehat{\bf E}_m|^2)^{1/2}
  \mathbb P(\Omega_{m}^{c})^{1/2}
  \leqslant
  \frac{\mathfrak c_2m^{1/2}}{\log(NT)}\mathbb P(\Omega_{m}^{c})^{1/2},
 \end{eqnarray*}
 where the constant $\mathfrak c_2 > 0$ doesn't depend on $m$ and $N$. Moreover,
 \begin{displaymath}
 \min_{\tau\in\mathcal S_m}\|\tau - b_I\|_{N}^{2}
 \leqslant\|b_I\|_{N}^{2}
 \quad {\rm because}\quad
 0\in\mathcal S_m,
 \end{displaymath}
 and then
 \begin{eqnarray*}
  \|\widehat b_m - b_I\|_{N}^{2}
  & = &
  \|\widehat b_m -\Pi_{N,m}(b_I)\|_{N}^{2} +
  \min_{\tau\in\mathcal S_m}\|\tau - b_I\|_{N}^{2}\\
  & \leqslant &
  \|\widehat b_m -\Pi_{N,m}(b_I)\|_{N}^{2} +\|b_I\|_{N}^{2}.
 \end{eqnarray*}
 Therefore,
 \begin{eqnarray*}
  \mathbb E(U_3) & \leqslant &
  \mathbb E(\|\widehat b_m -\Pi_{N,m}(b_I)\|_{N}^{2}
  \mathbf 1_{\Lambda_m\cap\Omega_{m}^{c}}) +
  \mathbb E(\|b_I\|_{N}^{2}\mathbf 1_{\Lambda_m\cap\Omega_{m}^{c}})\\
  & \leqslant &
  \frac{\mathfrak c_2m^{1/2}}{\log(NT)}
  \mathbb P(\Omega_{m}^{c})^{1/2} +
  \mathfrak c_1\mathbb P(\Omega_{m}^{c})^{1/2}.
 \end{eqnarray*}
\end{itemize}
So,
\begin{eqnarray*}
 \mathbb E(\|\widetilde b_m - b_I\|_{N}^{2})
 & \leqslant &
 \min_{\tau\in\mathcal S_m}\|b_I -\tau\|_{f_T}^{2}\\
 & &
 \hspace{1cm} +
 \frac{2m}{NT}(\|\sigma\|_{\infty}^{2} +
 \lambda\mathfrak c_{\zeta^2}\|\gamma\|_{\infty}^{2}) +
 \mathfrak c_2
 \frac{\sqrt{m\mathbb P(\Omega_{m}^{c})}}{\log(NT)} +
 \mathfrak c_1(
 \mathbb P(\Lambda_{m}^{c})^{1/2} +\mathbb P(\Omega_{m}^{c})^{1/2}).
\end{eqnarray*}
Therefore, by Lemma \ref{deviation_probabilities}, there exists a constant $\mathfrak c_3 > 0$, not depending on $m$ and $N$, such that
\begin{displaymath}
\mathbb E(\|\widetilde b_m - b_I\|_{N}^{2})
\leqslant
\min_{\tau\in\mathcal S_m}\|b_I -\tau\|_{f_T}^{2} +
\mathfrak c_3\left(\frac{m}{NT} +\frac{1}{N}\right).
\end{displaymath}
\hfill $\Box$
%

% Subsubsection : Proof of Lemma bound_variance_remainder.

%
\subsubsection{Proof of Lemma \ref{bound_variance_remainder}}\label{proof_variance_remainder}
In the sequel, the quadratic variation of any piecewise continuous stochastic process $(\Gamma_t)_{t\in [0,T]}$ is denoted by $(\llbracket\Gamma\rrbracket_t)_{t\in [0,T]}$. First of all, note that since $B$ and $Z$ are independent, for every $j\in\{1,\dots,m\}$,
\begin{eqnarray*}
 \left\llbracket\int_{0}^{.}
 \varphi_j(X_s)(\sigma(X_s)dB_s +\gamma(X_s)d\mathfrak Z_s)
 \right\rrbracket_T & = &
 \int_{0}^{T}\varphi_j(X_s)^2\sigma(X_s)^2ds +
 \int_{0}^{T}\varphi_j(X_s)^2\gamma(X_s)^2dZ_{s}^{(2)}\\
 & = &
 \int_{0}^{T}\varphi_j(X_s)^2(\sigma(X_s)^2 +
 \mathfrak c_{\zeta^2}\lambda\gamma(X_s)^2)ds\\
 & &
 \hspace{4cm} +
 \int_{0}^{T}\varphi_j(X_s)^2\gamma(X_s)^2d\mathfrak Z_{s}^{(2)}
\end{eqnarray*}
where, for every $t\in [0,T]$,
\begin{displaymath}
Z_{t}^{(2)} :=\llbracket\mathfrak Z\rrbracket_t =
\sum_{n = 1}^{\nu_t}\zeta_{n}^{2}
\quad {\rm and}\quad
\mathfrak Z_{t}^{(2)} := Z_{t}^{(2)} -\mathfrak c_{\zeta^2}\lambda t.
\end{displaymath}
By Jensen's inequality and Burkholder-Davis-Gundy's inequality (see Dellacherie et Meyer (1980), p. 303), there exists a constant $\mathfrak c_1 > 0$, not depending on $m$ and $N$, such that
\begin{eqnarray*}
 \mathbb E(|\widehat{\bf E}_{m}^{*}\widehat{\bf E}_m|^2) & \leqslant &
 m\sum_{j = 1}^{m}\mathbb E([\widehat{\bf E}_m]_{j}^{4})
 \leqslant
 \frac{\mathfrak c_1m}{N^4T^4}\sum_{j = 1}^{m}\mathbb E\left(
 \left\llbracket\sum_{i = 1}^{N}\int_{0}^{.}
 \varphi_j(X_{s}^{i})(\sigma(X_{s}^{i})dB_{s}^{i} +\gamma(X_{s}^{i})d\mathfrak Z_{s}^{i})
 \right\rrbracket_{T}^{2}\right)\\
 & \leqslant &
 \frac{2\mathfrak c_1m}{N^2T^4}\sum_{j = 1}^{m}\mathbb E\left(
 \left\llbracket\int_{0}^{.}
 \varphi_j(X_s)(\sigma(X_s)dB_s +\gamma(X_s)d\mathfrak Z_s)
 \right\rrbracket_{T}^{2}\right)\\
 & \leqslant &
 \frac{4\mathfrak c_1m}{N^2T^4}\sum_{j = 1}^{m}
 \left[\mathbb E\left[\left(\int_{0}^{T}\varphi_j(X_s)^2(\sigma(X_s)^2 +
 \mathfrak c_{\zeta^2}\lambda\gamma(X_s)^2)ds\right)^2\right]\right.\\
 & &
 \hspace{5cm}
 \left. +
 \mathbb E\left[\left(\int_{0}^{T}\varphi_j(X_s)^2
 \gamma(X_s)^2d\mathfrak Z_{s}^{(2)}\right)^2\right]
 \right].
\end{eqnarray*}
By Jensen's inequality,
\begin{eqnarray*}
 \left(\int_{0}^{T}\varphi_j(X_s)^2(\sigma(X_s)^2 +
 \mathfrak c_{\zeta^2}\lambda\gamma(X_s)^2)ds\right)^2
 & = &
 T^2\left(\int_{0}^{T}\varphi_j(X_s)^2(\sigma(X_s)^2 +
 \mathfrak c_{\zeta^2}\lambda\gamma(X_s)^2)\frac{ds}{T}\right)^2\\
 & \leqslant &
 T\int_{0}^{T}\varphi_j(X_s)^4(\sigma(X_s)^2 +
 \mathfrak c_{\zeta^2}\lambda\gamma(X_s)^2)^2ds.
\end{eqnarray*}
Moreover, since $\|\mathbf x\|_{4,m}\leqslant\|\mathbf x\|_{2,m}$ for every $\mathbf x\in\mathbb R^d$,
\begin{displaymath}
\sup_{x\in I}\sum_{j = 1}^{m}\varphi_j(x)^4 =
\sup_{x\in I}\|(\varphi_j(x))_j\|_{4,m}^{4}\leqslant
\sup_{x\in I}\|(\varphi_j(x))_j\|_{2,m}^{4}\leqslant L(m)^2.
\end{displaymath}
So, by applying twice the Fubini-Tonelli theorem,
\begin{eqnarray*}
 & &
 \sum_{j = 1}^{m}
 \mathbb E\left[\left(\int_{0}^{T}\varphi_j(X_s)^2(\sigma(X_s)^2 +
 \mathfrak c_{\zeta^2}\lambda\gamma(X_s)^2)ds\right)^2\right]\\
 & &
 \hspace{4cm}\leqslant
 2T\sum_{j = 1}^{m}\int_{0}^{T}
 \mathbb E(\varphi_j(X_s)^4(\sigma(X_s)^4 +
 \mathfrak c_{\zeta^2}^{2}\lambda^2\gamma(X_s)^4))ds\\
 & &
 \hspace{4cm} =
 2T^2\int_{-\infty}^{\infty}
 \underbrace{\left(\sum_{j = 1}^{m}\varphi_j(x)^4\right)}_{\leqslant L(m)^2}
 (\sigma(x)^4 +\mathfrak c_{\zeta^2}^{2}\lambda^2\gamma(x)^4)f_T(x)dx,
\end{eqnarray*}
and by the isometry type property of the stochastic integral with respect to $\mathfrak Z^{(2)}$,
\begin{eqnarray*}
 \sum_{j = 1}^{m}\mathbb E\left[\left(\int_{0}^{T}\varphi_j(X_s)^2
 \gamma(X_s)^2d\mathfrak Z_{s}^{(2)}\right)^2\right]
 & = &
 \lambda\mathfrak c_{\zeta^4}
 \sum_{j = 1}^{m}\int_{0}^{T}\mathbb E(\varphi_j(X_s)^4\gamma(X_s)^4)ds\\
 & \leqslant &
 \lambda\mathfrak c_{\zeta^4}TL(m)^2
 \int_{-\infty}^{\infty}\gamma(x)^4f_T(x)dx.
\end{eqnarray*}
Therefore,
\begin{displaymath}
\mathbb E(|\widehat{\bf E}_{m}^{*}\widehat{\bf E}_m|^2)
\leqslant
\frac{\mathfrak c_2}{N^2T^2}mL(m)^2
\end{displaymath}
with
\begin{displaymath}
\mathfrak c_2 =
8\mathfrak c_1\left(
\int_{-\infty}^{\infty}(\sigma(x)^4 +
\mathfrak c_{\zeta^2}^{2}\lambda^2\gamma(x)^4)f_T(x)dx +
\lambda\mathfrak c_{\zeta^4}\int_{-\infty}^{\infty}\gamma(x)^4f_T(x)dx
\right).
\end{displaymath}
\hfill $\Box$
%

% Subsubsection : Proof of Theorem th_model_selection.

%
\subsection{Proof of Theorem \ref{th_model_selection}} 
Let us consider the events
\begin{displaymath}
\Omega_N :=\bigcap_{m\in\mathcal M_{N}^{+}}\Omega_m
\quad {\rm and}\quad
\Xi_N :=\{\mathcal M_N\subset\widehat{\mathcal M}_N\subset\mathcal M_{N}^{+}\},
\end{displaymath}
where 
\begin{displaymath}
\mathcal M_{N}^{+} :=\left\{m\in\{1,\dots,N_T\} :
\mathfrak c_{\varphi}^{2}m(\|\mathbf\Psi_m^{-1}\|_{\rm op}^{2}\vee 1)
\leqslant 4\mathfrak d_T\frac{NT}{\log(NT)}\right\},
\end{displaymath}
and let us recall that
\begin{displaymath}
\mathfrak d_T =
\min\left\{\frac{\mathfrak c_T}{2},
\frac{1}{64\mathfrak c_{\varphi}^{2}T(\|f_T\|_{\infty} +\sqrt{\mathfrak c_T/2}/(3\mathfrak c_{\varphi}))}
\right\}.
\end{displaymath}
As a reminder, the sets $\widehat{\mathcal M}_N$ and $\mathcal M_N$ introduced in Section \ref{section_model_selection} are respectively defined by
\begin{displaymath}
\widehat{\mathcal M}_N =\left\{m\in\{1,\dots,N_T\} :
\mathfrak c_{\varphi}^{2}m(\|\widehat{\bf\Psi}_{m}^{-1}\|_{\rm op}^{2}\vee 1)
\leqslant\mathfrak d_T\frac{NT}{\log(NT)}\right\}
\end{displaymath}
and 
\begin{displaymath}
\mathcal M_N =\left\{m\in\{1,\dots,N_T\} :
\mathfrak c_{\varphi}^{2}m(\|\mathbf\Psi_{m}^{-1}\|_{\rm op}^{2}\vee 1)
\leqslant\frac{\mathfrak d_T}{4}\cdot\frac{NT}{\log(NT)}\right\}.
\end{displaymath}
The proof  of Theorem \ref{th_model_selection} relies on the three following lemmas.
%

% Lemma : Deviation probability for model selection.

%
\begin{lemma}\label{deviation_probability_model_selection}
Under Assumptions \ref{condition_m}, \ref{nested_models} and \ref{sub_exponential_Levy_measure}, there exists a constant $\mathfrak c_{\ref{deviation_probability_model_selection}} > 0$, not depending on $N$, such that
\begin{displaymath}
\mathbb P(\Xi_{N}^{c})\leqslant
\frac{\mathfrak c_{\ref{deviation_probability_model_selection}}}{N^6}.
\end{displaymath}
\end{lemma}
%

% Lemma : Bernstein type inequality.

%
\begin{lemma}\label{Bernstein_type_inequality}[Bernstein type inequality]
Consider the empirical process
\begin{displaymath}
\nu_N(\tau) :=
\frac{1}{NT}\sum_{i = 1}^{N}\int_{0}^{T}\tau(X_{s}^{i})(
\sigma(X_{s}^{i})dB_{s}^{i} +\gamma(X_{s}^{i})d\mathfrak Z_{s}^{i})
\textrm{ $;$ }
\tau\in\mathcal S_1\cup\dots\cup\mathcal S_{N_T}.
\end{displaymath}
Under Assumption \ref{sub_exponential_Levy_measure}, for every $\xi,v > 0$,
\begin{displaymath}
\mathbb P(\nu_N(\tau)\geqslant\xi,\|\tau\|_{N}^{2}\leqslant v^2)
\leqslant\exp\left[-\frac{NT\xi^2}{4[\mathfrak c_{\ref{Bernstein_type_inequality}}
(\|\sigma\|_{\infty}^{2} +
\|\gamma\|_{\infty}^{2})v^2 +\xi\|\tau\|_{\infty,I}\|\gamma\|_{\infty}]}\right]
\end{displaymath}
with
\begin{displaymath}
\mathfrak c_{\ref{Bernstein_type_inequality}} =
\frac{1}{2}\left[1\vee\int_{-\infty}^{\infty}e^{\mathfrak b|z|/2}\pi_{\lambda}(dz)\right].
\end{displaymath}
\end{lemma}
%

% Lemma : Bound on the empirical process.

%
\begin{lemma}\label{bound_empirical_process}
Under Assumptions \ref{condition_m} and \ref{sub_exponential_Levy_measure}, there exists a constant $\mathfrak c_{\ref{bound_empirical_process}} > 0$, not depending on $N$, such that for every $m\in\mathcal M_N$,
\begin{displaymath}
\mathbb E\left[\left(\left[\sup_{\tau\in\mathcal B_{m,m'}}\nu_n(\tau)\right]^2 -
p(m,\widehat m)\right)_{+}\mathbf 1_{\Xi_N\cap\Omega_N}\right]\leqslant
\frac{\mathfrak c_{\ref{bound_empirical_process}}}{NT}
\end{displaymath}
where, for every $m'\in\mathcal M_N$,
\begin{displaymath}
\mathcal B_{m,m'} :=
\left\{\tau\in\mathcal S_{m\wedge m'} :\|\tau\|_{f_T} = 1\right\}
\quad {\rm and}\quad
p(m,m') :=\frac{\mathfrak c_{\mathrm{cal}}}{8}\cdot\frac{m\vee m'}{NT}.
\end{displaymath}
\end{lemma}
The proof of Lemma \ref{Bernstein_type_inequality} is postponed to Subsubsection \ref{proof_Bernstein_type_inequality}. Lemma \ref{bound_empirical_process} is a consequence of Lemma \ref{Bernstein_type_inequality} thanks to the $\mathbb L_{f_T}^{2}$-$\mathbb L^{\infty}$ chaining technique (see Comte (2001), Proposition 4). Finally, the proof of Lemma \ref{deviation_probability_model_selection} remains the same as the proof of Comte and Genon-Catalot (2020b), Eq. (6.17), because $(B^1,\mathfrak Z^1),\dots,(B^N,\mathfrak Z^N)$ are independent.
%

% Subsubsection : Steps of the proof.

%
\subsubsection{Steps of the proof}
First of all,
\begin{eqnarray}
 \|\widehat b_{\widehat m} - b_I\|_{N}^{2} & = &
 \|\widehat b_{\widehat m} - b_I\|_{N}^{2}\mathbf 1_{\Xi_{N}^{c}} +
 \|\widehat b_{\widehat m} - b_I\|_{N}^{2}\mathbf 1_{\Xi_N}
 \nonumber\\
 \label{bound_adaptive_estimator_1}
 & =: &
 U_1 + U_2.
\end{eqnarray}
Let us find suitable bounds on $\mathbb E(U_1)$ and $\mathbb E(U_2)$.
\begin{itemize}
 \item {\bf Bound on $\mathbb E(U_1)$.} Since
 \begin{displaymath}
 \|\widehat b_{\widehat m} -\Pi_{N,\widehat m}(b_I)\|_{N}^{2} =
 \widehat{\bf E}_{\widehat m}^{*}\widehat{\bf\Psi}_{\widehat m}^{-1}
 \widehat{\bf E}_{\widehat m},
 \end{displaymath}
 by the definition of $\widehat{\mathcal M}_N$, and by Lemma \ref{bound_variance_remainder},
 \begin{eqnarray*}
  \mathbb E(\|\widehat b_{\widehat m} -\Pi_{N,\widehat m}(b_I)\|_{N}^{2}
  \mathbf 1_{\Xi_{N}^{c}})
  & \leqslant &
  \mathbb E(\|\widehat{\bf\Psi}_{\widehat m}^{-1}\|_{\rm op}
  |\widehat{\bf E}_{NT}^{*}\widehat{\bf E}_{NT}|
  \mathbf 1_{\Xi_{N}^{c}})\\
  & \leqslant &
  \left[
  \mathfrak d_T
  \frac{NT}{\log(NT)}\right]^{1/2}
  \mathbb E(|\widehat{\bf E}_{NT}^{*}\widehat{\bf E}_{NT}|^2)^{1/2}
  \mathbb P(\Xi_{N}^{c})^{1/2}
  \leqslant
  \frac{\mathfrak c_1N}{\log(NT)}\mathbb P(\Xi_{N}^{c})^{1/2},
 \end{eqnarray*}
 where the constant $\mathfrak c_1 > 0$ doesn't depend on $N$. Then,
 \begin{eqnarray*}
  \mathbb E(U_1) & \leqslant &
  \mathbb E(\|\widehat b_{\widehat m} -\Pi_{N,\widehat m}(b_I)\|_{N}^{2}
  \mathbf 1_{\Xi_{N}^{c}}) +
  \mathbb E(\|b_I\|_{N}^{2}\mathbf 1_{\Xi_{N}^{c}})\\
  & \leqslant &
  \frac{\mathfrak c_1N}{\log(NT)}\mathbb P(\Xi_{N}^{c})^{1/2} +
  \mathfrak c_2\mathbb P(\Xi_{N}^{c})^{1/2}
 \end{eqnarray*}
 with
 \begin{displaymath}
 \mathfrak c_2 =\left(\int_{-\infty}^{\infty}b_I(x)^4f_T(x)dx\right)^{1/2}.
 \end{displaymath}
 So, by Lemma \ref{deviation_probability_model_selection}, there exists a constant $\mathfrak c_3 > 0$, not depending on $N$, such that
 \begin{displaymath}
 \mathbb E(U_1)
 \leqslant\frac{\mathfrak c_3}{N}.
 \end{displaymath}
 \item {\bf Bound on $\mathbb E(U_2)$.} Note that
 \begin{eqnarray*}
  U_2 & = &
  \|\widehat b_{\widehat m} - b_I\|_{N}^{2}\mathbf 1_{\Xi_N\cap\Omega_{N}^{c}} +
  \|\widehat b_{\widehat m} - b_I\|_{N}^{2}\mathbf 1_{\Xi_N\cap\Omega_N}\\
  & =: &
  U_{2,1} + U_{2,2}.
 \end{eqnarray*}
 On the one hand, by Lemma \ref{deviation_probabilities}, there exists a constant $\mathfrak c_4 > 0$, not depending on $N$, such that
 \begin{displaymath}
 \mathbb P(\Xi_N\cap\Omega_{N}^{c})
 \leqslant
 \sum_{m\in\mathcal M_{N}^{+}}\mathbb P(\Omega_{m}^{c})
 \leqslant\frac{\mathfrak c_4}{N^6}.
 \end{displaymath}
 Then, as for $\mathbb E(U_1)$, there exists a constant $\mathfrak c_5 > 0$, not depending on $N$, such that
 \begin{displaymath}
 \mathbb E(U_{2,1})
 \leqslant\frac{\mathfrak c_5}{N}.
 \end{displaymath}
 On the other hand,
 \begin{displaymath}
 \gamma_N(\tau') -\gamma_N(\tau) =
 \|\tau' - b\|_{N}^{2} -\|\tau - b\|_{N}^{2} - 2\nu_N(\tau' -\tau)
 \end{displaymath}
 for every $\tau,\tau'\in\mathcal S_1\cup\dots\cup\mathcal S_{N_T}$. Moreover, since
 \begin{displaymath}
 \widehat m =
 \arg\min_{m\in\widehat{\mathcal M}_N}\{-\|\widehat b_m\|_{N}^{2} +
 {\rm pen}(m)\} =
 \arg\min_{m\in\widehat{\mathcal M}_N}\{\gamma_N(\widehat b_m) +
 {\rm pen}(m)\},
 \end{displaymath}
 for every $m\in\widehat{\mathcal M}_N$,
 \begin{equation}\label{bound_adaptive_estimator_2}
 \gamma_N(\widehat b_{\widehat m}) + {\rm pen}(\widehat m)
 \leqslant
 \gamma_N(\widehat b_m) + {\rm pen}(m).
 \end{equation}
 On the event $\Xi_N =\{\mathcal M_N\subset\widehat{\mathcal M}_N\subset\mathcal M_{N}^{+}\}$, Inequality (\ref{bound_adaptive_estimator_2}) remains true for every $m\in\mathcal M_N$. Then, on $\Xi_N$, for any $m\in\mathcal M_N$, since $\mathcal S_m +\mathcal S_{\widehat m}\subset\mathcal S_{m\vee\widehat m}$ under Assumption \ref{nested_models},
 \begin{eqnarray*}
  \|\widehat b_{\widehat m} - b_I\|_{N}^{2}
  & \leqslant &
  \|\widehat b_m - b_I\|_{N}^{2} +
  2\|\widehat b_{\widehat m} -\widehat b_m\|_{f_T}
  \nu_N\left(\frac{\widehat b_{\widehat m} -\widehat b_m}{
  \|\widehat b_{\widehat m} -\widehat b_m\|_{f_T}}\right) +
  {\rm pen}(m) - {\rm pen}(\widehat m)\\
  & \leqslant &
  \|\widehat b_m - b_I\|_{N}^{2} +
  \frac{1}{8}\|\widehat b_{\widehat m} -\widehat b_m\|_{f_T}^{2}\\
  & &
  \hspace{1.5cm} +
  8\left(\left[\sup_{\tau\in\mathcal B_{m,\widehat m}}|\nu_N(\tau)|\right]^2
  - p(m,\widehat m)\right)_+ +
  {\rm pen}(m) + 8p(m,\widehat m) - {\rm pen}(\widehat m).
 \end{eqnarray*}
 Since $\|.\|_{f_T}^{2}\mathbf 1_{\Omega_N}\leqslant 2\|.\|_{N}^{2}\mathbf 1_{\Omega_N}$ on $\mathcal S_1\cup\dots\cup\mathcal S_{\max(\mathcal M_{N}^{+})}$, and since $8p(m,\widehat m)\leqslant {\rm pen}(m) + {\rm pen}(\widehat m)$, on $\Xi_N\cap\Omega_N$,
 \begin{displaymath}
 \|\widehat b_{\widehat m} - b_I\|_{N}^{2}
 \leqslant
 3\|\widehat b_m - b_I\|_{N}^{2} + 4{\rm pen}(m) +
 16\left(\left[\sup_{\tau\in\mathcal B_{m,\widehat m}}|\nu_N(\tau)|\right]^2
 - p(m,\widehat m)\right)_+.
 \end{displaymath}
 So, by Lemma \ref{bound_empirical_process},
 \begin{eqnarray*}
  \mathbb E(U_{2,2})
  & \leqslant &
  \min_{m\in\mathcal M_N}
  \{\mathbb E(3\|\widehat b_m - b_I\|_{N}^{2}\mathbf 1_{\Xi_N}) + 4{\rm pen}(m)\} +
  \frac{16\mathfrak c_{\ref{bound_empirical_process}}}{NT}\\
  & \leqslant &
  \mathfrak c_6\min_{m\in\mathcal M_N}
  \left\{\inf_{\tau\in\mathcal S_m}
  \|\tau - b_I\|_{f_T}^{2} +\frac{m}{NT}
  \right\} +
  \frac{\mathfrak c_6}{N}
 \end{eqnarray*}
 where $\mathfrak c_6 > 0$ is a deterministic constant not depending on $N$.
\end{itemize}
\hfill $\Box$
%

% Subsubsection : Proof of Lemma Bernstein_type_inequality.

%
\subsubsection{Proof of Lemma \ref{Bernstein_type_inequality}}\label{proof_Bernstein_type_inequality}
Consider $\tau\in\mathcal S_1\cup\dots\cup\mathcal S_{N_T}$ and, for any $i\in\{1,\dots,N\}$, let $M^i(\tau) = (M^i(\tau)_t)_{t\in [0,T]}$ be the martingale defined by
\begin{displaymath}
M^i(\tau)_t :=
\int_{0}^{t}\tau(X_{s}^{i})(\sigma(X_{s}^{i})dB_{s}^{i} +
\gamma(X_{s}^{i})d\mathfrak Z_{s}^{i})
\textrm{ $;$ }
\forall t\in [0,T].
\end{displaymath}
Moreover, for every $\varepsilon > 0$, consider
\begin{displaymath}
Y_{\varepsilon}^{i}(\tau) :=
\varepsilon M^i(\tau) - A_{\varepsilon}^{i}(\tau) - B_{\varepsilon}^{i}(\tau),
\end{displaymath}
where $A_{\varepsilon}^{i}(\tau) = (A_{\varepsilon}^{i}(\tau)_t)_{t\in [0,T]}$ and $B_{\varepsilon}^{i}(\tau) = (B_{\varepsilon}^{i}(\tau)_t)_{t\in [0,T]}$ are the stochastic processes defined by
\begin{eqnarray*}
 A_{\varepsilon}^{i}(\tau)_t & := &
 \frac{\varepsilon^2}{2}\int_{0}^{t}\tau(X_{s}^{i})^2\sigma(X_{s}^{i})^2ds\\
 & &
 \hspace{1cm} {\rm and}\quad
 B_{\varepsilon}^{i}(\tau)_t :=
 \int_{0}^{t}\left[\int_{-\infty}^{\infty}(
 e^{\varepsilon z\tau(X_{s}^{i})\gamma(X_{s}^{i})} -
 \varepsilon z\tau(X_{s}^{i})\gamma(X_{s}^{i}) - 1)\pi_{\lambda}(dz)\right]ds
\end{eqnarray*}
for every $t\in [0,T]$. The proof is dissected in three steps.
\\
\\
{\bf Step 1.} Note that for any $i\in\{1,\dots,N\}$ and $t\in [0,T]$,
\begin{displaymath}
|\tau(X_{t}^{i})\gamma(X_{t}^{i})|\leqslant
\|\tau\|_{\infty,I}\|\gamma\|_{\infty}
\end{displaymath}
and then, by Assumption \ref{sub_exponential_Levy_measure},
\begin{displaymath}
\mathbb E\left(\int_{0}^{t}\int_{|z| > 1}
|e^{\varepsilon z\tau(X_{s}^{i})\gamma(X_{s}^{i})} - 1|\pi_{\lambda}(dz)ds\right) <\infty
\end{displaymath}
for any $\varepsilon\in (0,\varepsilon^*)$ with $\varepsilon^* = (\mathfrak b\wedge 1)/(2\|\tau\|_{\infty,I}\|\gamma\|_{\infty})$. So, $(\exp(Y_{\varepsilon}^{i}(\tau)_t))_{t\in[0,T]}$ is a local martingale by Applebaum (2009), Corollary 5.2.2. In other words, there exists an increasing sequence of stopping times $(T_{n}^{i})_{n\in\mathbb N}$ such that $\lim_{n\rightarrow\infty}T_{n}^{i} =\infty$ a.s. and $(\exp(Y_{\varepsilon}^{i}(\tau)_{t\wedge T_{n}^{i}})_{t\in [0,T]}$ is a martingale. Therefore, by Lebesgue's theorem and Markov's inequality, for every $\rho > 0$, the stochastic process $Y_{N,\varepsilon}(\tau) := Y_{\varepsilon}^{1}(\tau) +\dots + Y_{\varepsilon}^{N}(\tau)$ satisfies
\begin{eqnarray*}
 \mathbb P(e^{Y_{N,\varepsilon}(\tau)_T} >\rho) & = &
 \lim_{n\rightarrow\infty}\mathbb P\left(\exp\left[
 \sum_{i = 1}^{N}Y_{\varepsilon}^{i}(\tau)_{T\wedge T_{n}^{i}}\right] >\rho\right)\\
 & \leqslant &
 \frac{1}{\rho}\lim_{n\rightarrow\infty}
 \mathbb E(\exp(Y_{\varepsilon}^{1}(\tau)_{T\wedge T_{n}^{1}}))^N =
 \frac{1}{\rho}\mathbb E(\exp(Y_{\varepsilon}^{1}(\tau)_0))^N =\frac{1}{\rho}.
\end{eqnarray*}
{\bf Step 2.} For any $\varepsilon\in (0,\varepsilon^*)$ and $t\in [0,T]$, let us find suitable bounds on
\begin{displaymath}
A_{N,\varepsilon}(\tau)_t :=\sum_{i = 1}^{N}A_{\varepsilon}^{i}(\tau)_t
\quad {\rm and}\quad
B_{N,\varepsilon}(\tau)_t :=\sum_{i = 1}^{N}B_{\varepsilon}^{i}(\tau)_t.
\end{displaymath}
On the one hand,
\begin{equation}\label{Bernstein_type_inequality_1}
A_{N,\varepsilon}(\tau)_t\leqslant
\frac{\varepsilon^2\|\sigma\|_{\infty}^{2}}{2}
\sum_{i = 1}^{N}\int_{0}^{t}\tau(X_{s}^{i})^2ds
\leqslant\frac{\varepsilon^2\|\sigma\|_{\infty}^{2}\|\tau\|_{N}^{2}NT}{2}.
\end{equation}
On the other hand, for every $\beta\in (-\mathfrak b/2,\mathfrak b/2)$, by Taylor's formula and Assumption \ref{sub_exponential_Levy_measure},
\begin{eqnarray*}
 \int_{-\infty}^{\infty}(e^{\beta z} -\beta z - 1)\pi_{\lambda}(dz) & = &
 \beta^2\int_{-\infty}^{\infty}
 \left(\int_{0}^{1}(1 -\theta)e^{\theta\beta z}d\theta\right)\pi_{\lambda}(dz)\\
 & \leqslant &
 \frac{\mathfrak c_1}{2}\beta^2
 \quad {\rm with}\quad
 \mathfrak c_1 =
 \int_{-\infty}^{\infty}e^{\mathfrak b|z|/2}\pi_{\lambda}(dz) <\infty.
\end{eqnarray*}
Since $\varepsilon\in (0,\varepsilon^*)$, one can take $\beta =\varepsilon\tau(X_{s}^{i})\gamma(X_{s}^{i})$ for any $s\in [0,t]$ and $i\in\{1,\dots,N\}$, and then
\begin{equation}\label{Bernstein_type_inequality_2}
B_{N,\varepsilon}(\tau)_t\leqslant
\frac{\mathfrak c_1\varepsilon^2}{2}\sum_{i = 1}^{N}\int_{0}^{t}
\tau(X_{s}^{i})^2\gamma(X_{s}^{i})^2ds
\leqslant\frac{\mathfrak c_1\varepsilon^{2}\|\gamma\|_{\infty}^{2}\|\tau\|_{N}^{2}NT}{2}.
\end{equation}
Therefore, Inequalities (\ref{Bernstein_type_inequality_1}) and (\ref{Bernstein_type_inequality_2}) lead to
\begin{displaymath}
A_{N,\varepsilon}(\tau)_t + B_{N,\varepsilon}(\tau)_t\leqslant
\mathfrak c_2\varepsilon^2(\|\sigma\|_{\infty}^{2} +
\|\gamma\|_{\infty}^{2})\|\tau\|_{N}^{2}NT
\quad {\rm with}\quad
\mathfrak c_2 =
\frac{1}{2}(1\vee\mathfrak c_1).
\end{displaymath}
{\bf Step 3 (conclusion).} Consider $M_N(\tau) := M^1(\tau) +\dots + M^N(\tau)$. For any $\varepsilon\in (0,\varepsilon^*)$ and $\xi,v > 0$, thanks to Step 2,
\begin{eqnarray*}
 \mathbb P(\nu_N(\tau)\geqslant\xi,\|\tau\|_{N}^{2}\leqslant v^2)
 & \leqslant &
 \mathbb P(e^{\varepsilon M_N(\tau)_T}\geqslant e^{NT\varepsilon\xi},
 A_{N,\varepsilon}(\tau)_T + B_{N,\varepsilon}(\tau)_T\leqslant
 \mathfrak c_2\varepsilon^2(\|\sigma\|_{\infty}^{2} +\|\gamma\|_{\infty}^{2})NTv^2)\\
 & \leqslant &
 \mathbb P(e^{Y_{N,\varepsilon}(\tau)_T}\geqslant
 \exp(NT\varepsilon\xi -\mathfrak c_2\varepsilon^2(\|\sigma\|_{\infty}^{2} +\|\gamma\|_{\infty}^{2})NTv^2)).
\end{eqnarray*}
Moreover, taking
\begin{displaymath}
\varepsilon =
\frac{\xi}{2\mathfrak c_2(\|\sigma\|_{\infty}^{2} +\|\gamma\|_{\infty}^{2})v^2 +\xi/\varepsilon^*}
<\varepsilon^*
\end{displaymath}
leads to
\begin{eqnarray*}
 NT\varepsilon\xi -\mathfrak c_2\varepsilon^2(\|\sigma\|_{\infty}^{2} +
 \|\gamma\|_{\infty}^{2})NTv^2
 & = &
 \frac{NT\xi^2[\mathfrak c_2(\|\sigma\|_{\infty}^{2} +\|\gamma\|_{\infty}^{2})v^2 +
 \xi/\varepsilon^*]}{[2\mathfrak c_2(\|\sigma\|_{\infty}^{2} +\|\gamma\|_{\infty}^{2})v^2 +
 \xi/\varepsilon^*]^2}\\
 & \geqslant &
 \frac{NT\xi^2}{4[\mathfrak c_2(\|\sigma\|_{\infty}^{2} +
 \|\gamma\|_{\infty}^{2})v^2 +\xi/\varepsilon^*]}.
\end{eqnarray*}
Therefore, by Step 1,
\begin{displaymath}
\mathbb P(\nu_N(\tau)\geqslant\xi,\|\tau\|_{N}^{2}\leqslant v^2)
\leqslant
\exp\left(-\frac{NT\xi^2}{4[\mathfrak c_2(\|\sigma\|_{\infty}^{2} +
\|\gamma\|_{\infty}^{2})v^2 +\xi\|\tau\|_{\infty,I}\|\gamma\|_{\infty}]}\right).
\end{displaymath}
\hfill $\Box$
%

% Section : Figures and tables.

%
\section{Figures and tables}\label{section_tables_figures}
\begin{table}[!h]
\begin{center}
\begin{tabular}{|l||c|c|c|}
 \hline
  & Model 1 & Model 2 & Model 3\\
 \hline
 \hline
 Mean MISE & 0.1251 & 0.1469 & 0.1825\\
 \hline
 StD MISE & 0.0950 & 0.1688 & 0.1928\\
 \hline
\end{tabular}
\medskip
\caption{Means and StD of the MISE of $\widehat b_{\widehat m}$ (100 repetitions).}\label{table_MISE}
\end{center}
\end{table}
\newpage
%

% References.

%

%
\end{document}